# Revisit the scheduling problem in Hurdle, V.F., 1973: A novel mathematical solution approach and two extensions

Wenbo Fan

School of Transportation and Logistics, National Engineering Laboratory of Integrated Transportation Big Data Application Technology, Southwest Jiaotong University, Chengdu, China; Email: wbfan@swjtu.edu.cn

**Abstract**

The scheduling problem in Hurdle (1973) was formulated in a general form that simultaneously concerned the vehicle dispatching, circulating, fleet sizing, and patron queueing. As a constrained variational problem, it remains not fully solved for decades. With technical prowess in graphic analysis, the author unveiled the closed-form solution for the optimal dispatch rates (with key variables undetermined though), but only suggested the lower and upper bounds of the optimal fleet size. Additionally, such a graphic analysis method lacks high efficiency in computing specific scheduling problems, which are often of a large scale (e.g., hundreds of bus lines). In light of this, the paper proposes a novel mathematical solution approach that first relaxes the original problem to an unconstrained one, and then attacks it using calculus of variations. The corresponding Euler-Lagrange equation yields the closed-form solution of the optimal dispatching rates, to which Hurdle's results are a special case. Thanks to the proposed approach, the optimal fleet size can also be solved. This paper completes the work of Hurdle (1973) by formalizing a solution method and generalizing the results. Based on that, we further make two extensions to the scheduling problem of general bus lines with multiple origins and destinations and that of mixed-size or modular buses. Closed-form results are also obtained with new insights. Among others, we find that the solutions for shuttle/feeder lines are a special case to our results of general bus lines. Numerical examples are also provided to demonstrate the effectiveness and efficiency of the proposed approach.

## 1. Introduction

Scheduling problem is one of the fundamental problems of scheduled transportation systems such as buses, subways, railroads, airlines, shipping lines, and the postal system. The problem is "simple" by description: to design the time schedule of bus/train/airplane/carrier dispatches for a certain level of service (LOS). Its mathematical problems, however, could be very challenging to resolve. Additionally, scheduling problems are often associated with the fleet sizing subproblem, which seeks for the minimum number of vehicles that carry out the schedule plan and circulate in service. To concentrate, we focus our following discussion on the scheduling problems of transit, e.g., buses and subways, and start from a literature review before diving into details. Such an overview allows us to better understand the big picture to which the work of Hurdle (1973) belongs and the necessity of our revisiting.

## 1.1. Literature review

The literature of transit scheduling studies can be generally categorized into two groups that consider the time-invariant demand and time-varying demand, respectively. The former is often set for planning tasks, e.g., transit network designs (Ceder and Wilson, 1986; Ibarra-Rojas et al., 2015), where the average peak-hour demand is chosen to be served and the scheduling problem is embedded as a subproblem (also named frequency setting problem). The subproblem involves a limited number of continuous decision variables, e.g., frequencies or headways of lines in the network, and can be solved according to the vehicle capacity or desired vehicle load constraints (Kepaptsoglou and Karlaftis, 2009). The study challenge or difficulty often lies in the integration with other planning-related problems, e.g., single-mode and multimodal network designs (Arbex and da Cunha, 2015; Verbas, Mahmassani, 2015; Szeto and Jiang, 2014), passenger assignment (Ibarra-Rojas et al., 2019), scheduling coordination/synchronization (Liu and Ceder, 2018), and the allocation of modular vehicles (Pei et al., 2021); see Guihaire and Hao (2008) for a good survey.

The second group of studies targets at daily tactical operations that concern the temporal variation of demand. The corresponding scheduling problems become very challenging due to the complexity in modeling the cost metrics (e.g., waiting times) of the time-varying and occasionally oversaturated demand. When oversaturated, some patrons inevitably have to wait in queue for more than one arrival of vehicles to get boarded/picked up. This queueing phenomenon interweaves the decisions of dispatches in distant times (Chen et al., 2020). The present study belongs to this group. Thus, a more detailed survey is presented below regarding the works that considered time-varying and/or queueing patrons.

According to the research methodologies, the time-varying scheduling studies can be classified into the ones of Continuum Approximation (CA) method and the others of mathematical modeling approach. The former models the problem using a continuous dispatching rate/frequency (or headway) as a function of time; while the latter models the dispatch times (i.e., clock times) of each vehicle in the fleet as decision variables, which can be modeled as continuous variables or binary variables indicating a dispatch or not at that time.

In the first class, Newell (1971) firstly formulated an analytical scheduling model for a shuttle/feeder bus line to minimize the total patrons waiting time. The author derived both an exact mathematical model (of which the decision variables were the dispatch time of each vehicle) and a CA model with the dispatching rate being a single decision function of time. The closed-form solution was found, which insightfully unveiled that the optimal dispatching rate should be proportional to the square root of patrons' arrival rate. This seminal work contained several assumptions that ignored the circulation of buses and the queueing of patrons. Extensions had been made by a number of follow-up studies, as summarized in Table 1. Considering a round-trip bus line, Salzborn (1972) firstly determined the fleet size by the peak-period demand, and then obtained the optimal dispatching rate by solving the fleet size constrained waiting time minimization problem via calculus of variations. Next, a major extension was made by Hurdle (1973), who formulated the problem as the minimization of the total cost for patrons' waiting time and operator's fleet purchase and operation. The proposed model explicitly accounted for the fleet size constraint and patrons queueing dynamics (see Section 2 for details). Such a problem was complicated in the form of a constrained variational problem and attacked through graphic analysis. With technical prowess, the author optimized the dispatching rates and obtained lower and upper bounds for the optimal fleet size. It was suggested to be findable via a trial-and-error method.

**Table 1. Representative CA-based studies on transit scheduling problems.**

| Works | Major contributions | Limitations |
|---|---|---|
| Newell (1971) | • Closed-form solution for optimal dispatching rate | No fleet sizing |
| Salzborn (1972) | • Fleet size determined by the peak-period demand<br>• Closed-form solutions for peak-period and off-peak dispatching rates | No queues allowed |
| Hurdle (1973) | • Accounted for queues<br>• Closed-form solutions for peak-period and off-peak dispatching rates<br>• Lower, upper bounds of the optimal fleet size | Optimal dispatching rate not fully solved; no solution of the optimal fleet size; |
| Chen et al. (2020) | • Optimal dispatching policy for modular vehicles | No fleet sizing |
| **This paper** | • Calculus of variations used to find the closed-form solutions of the optimal peak-period and off-peak dispatching rates<br>• The optimal fleet size is solved<br>• Two extensions | |

Later on, the problem of scheduling for time-varying demand was mainly studied by using mathematical models partly due to the difficulty in solving CA models. Examples in this realm include Sheffi and Sugiyama (1982), Ceder (1984), Niu and Zhou (2013), Barrena et al. (2014a, b), Niu et al. (2015), Yin et al. (2017), and Bertsimas et al. (2020) among others. These studies overcame the limits of CA models by considering more realistic and sophisticated conditions. For instance, Sheffi and Sugiyama (1982) extended Newell's (1971) mathematical model to a general bus line with multiple origins and destinations, and further considered non-fixed bus travel time and stochastic arrivals of patrons. Ceder (1984) made use of surveyed passenger data and accordingly determined the dispatching rate for different periods so as to achieve a certain LOS (e.g., even headway or even vehicle load). Recent time-varying scheduling works include Yin et al (2017) that considered energy efficiency of trains, Shi et al. (2020) that considered common bus lines with modular vehicles, and Bertsimas, et al. (2020) that jointly optimized dispatching rates and pricing for lines in a multimodal transit network.

In the field of railway, scheduling problem is also in the name of timetabling problem that often emphasizes the vehicle operation rather than the demand-oriented service design; see Caprara et al (2002) for a good survey. Barrena et al. (2014a) can be seen as an extension to Sheffi and Sugiyama (1982) by optimizing the dispatch/departure times of vehicles/trains at not only the depot but also all stations along the rail line under fleet size constraint. The corresponding problem is NP-hard and solved via a metaheuristic algorithm (i.e., adaptive large neighborhood search). An equivalent formulation of the problem was also proposed in Barrena et al. (2014b), but limited to small instances. Niu et al. (2015) tackled a similar problem in a skip-stop rail corridor.

The first scheduling work that explicitly tackled oversaturated demand may be Niu and Zhou (2013) for an urban rail line subject to the train fleet size constraint. To capture the queueing phenomenon at stations, the authors divided the time horizon into very small intervals. Consequently, the problem size was significantly enlarged with each time point being the potential departure time of trains. A genetic algorithm was used to solve the binary integer nonlinear programming model with no guarantee of finding the global optimum. In the study of rescheduling a disrupted metro line under over-crowed situations, Gao et al. (2016) applied linearization techniques to deal with the queueing-related variables and formulated a mixed integer linear programming model, which was then solved via a heuristic iterative algorithm. Shang et al. (2019) turned attention to the equity-oriented scheduling for an oversaturated skip-stop rail corridor. They formulated passengers and train stopping patterns in a time-space-state multi-commodity flow-modeling framework, which introduced additional dimensions and numerous binary variables. The proposed linear programming model was decomposed into least-cost subproblems under Lagrangian relaxation and then solved by Bellman-Ford algorithm. Shi et al. (2018) jointly modeled the train timetabling problem with collaborative passenger flow control on an oversaturated metro line. Recently, a notable exception to the most mathematical modeling studies is Chen et al. (2020), who returned to CA method for the scheduling problem of a novel shuttle service with modular vehicles and oversaturated demand. The formulated problem ignored the fleet sizing issue, however.

From the above survey, it can be seen that although with application conditions substantially extended, the mathematical scheduling models often require enormous efforts in data preparation, model formulation, and solution computation. Consequently, they apply only to particular problems. Little new insight can be drawn from these sophisticated models. More importantly, a theoretical gap remains unfilled since first proposed in Hurdle (1973): the joint optimization of the dispatches and the fleet size, which witnesses no success in the above mathematical models. Without the optimal fleet size, the optimized dispatching rate or schedule times are suboptimal solutions.

Although seemingly slowly developed in the scheduling problem, the CA-based models have seen a rapid growth in the general field of transportation and logistics; see examples in transit network design including Daganzo (2010), Ouyang et al. (2014), Chen et al. (2015), Fan et al. (2018), etc.; and examples in facility location (Carlsson and Jia, 2014) and integrated supply chain and logistics (Carlsson and Song, 2017) among others. A recent review by Ansari et al. (2018) provides a more comprehensive survey on CA models. Therefore, we feel necessary and promising to revisit the underdeveloped CA-based scheduling problems (Hurdle's problem in particular) and seek opportunities for further extensions.

1.2. The importance and limitation of Hurdle's scheduling problem

We highlight the importance and limitation of Hurdle's scheduling problem by three points, which also explain why we pay the revisit.

    i.  Hurdle's scheduling problem is an open question that remains not fully solved for 48 years upon the time of writing. As said, the fleet sizing problem was not resolved. The suggested trial-and-error method would be time consuming to find the optimal fleet size. Additionally, the optimized dispatching rates also contained unknown variables (i.e., the queue beginning time and dissipating time) to be determined. The graphic analysis used in Hurdle (1973) is good for visualization and understanding but lacks high efficiency in solving specific large-scaled problems. More importantly, the optimized dispatching rates obtained via graphic analysis needs a rigorous proof. In

doing so, we find that Hurdle's result is a special case to the optimal solution derived by our mathematical solution approach (Section 3.1).

ii. Hurdle's scheduling problem is of theoretical importance even to this day. The analytically formulated problem holds the promise of closed-form solution, which would unveil fundamental insights into the cause-and-effect relation between the optimal design and key factors (including the demand profiles, patrons' values of time, capital and operational costs of buses, etc.). In addition, the modeling framework is in a general form that welcomes applications to broader and modern scheduling problems. For instance, the shuttle/feeder-service-oriented model can be extended to general bus lines with demand of multiple origins and destinations (see Section 4.1); the modeling of dispatches in terms of seats/places for patrons can be applied to mixed-size buses and futuristic modular buses according to the associated operation cost structure (see Section 4.2); and the consideration of buses' circulation/round-trip time may also include that (a layover time for charging) of electric buses, which have been promoted to tackle the transportation-related environmental problems of the time in many cities around the world.

iii. The solution to Hurdle's scheduling problem is also of practical usefulness. Specific times of dispatches can be readily furnished (via integration of the dispatching rate; see two recipes in Daganzo (2005) and Chen et al. (2020), respectively) for implementation. Fine tunes can be made to satisfy further practical concerns (such as 'clock headway'), without much deteriorating the system performance, because it is rather flat in the neighborhood of the optimum (Daganzo, 1997).

## 1.3 Contributions of this paper

By revisiting the scheduling problem in Hurdle (1973), this paper contributes to the literature in four-fold:

i. A novel mathematical solution approach is proposed to fully resolve Hurdle's (1973) constrained scheduling problem. Specifically, the peak-period dispatching rate is solved to the optimal with Hurdle's result being a special case. The optimal fleet size, that was not determined in Hurdle (1973), is now solved. Useful system metrics, e.g., the queueing profile, are also obtained with high efficiency.

ii. Thanks to the proposed approach, closed-form solutions are obtained for the optimal dispatching rates of general bus lines. Simultaneously, the optimal fleet size can also be solved. To our findings, the results for shuttle/feeder bus lines become a special case.

iii. The scheduling problem of mixed-size/modular buses is also analytically studied with closed-form solutions obtained for the optimal vehicle sizes.

iv. New insights are unravelled from the above analytical results and a variety of numerical experiments. Among others, we find that the optimal peak-period dispatch rate of general bus lines is oftentimes determined by the maximum cross-sectional demand flux at all stops.

Next, we restate Hurdle's scheduling problem and solution as the starting point for developing our mathematical solution approach (Section 3) and extensions (Section 4). Section 5 concludes the paper with discussions of study limitations and future research directions.

## 2. Hurdle's scheduling problem restatement

Consider a single bus line of shuttle or feeder services for demand of one-to-one, one-to-many, or many-to-one patterns. The one-to-one demand pattern indicates that all patrons board and alight at the beginning and ending stops. For one-to-many pattern, all patrons board at the beginning stop and alight at other stops along the line. In the case of many-to-one demand, patrons board at stops along the line and all alight at the last stop. For the above three demand patterns, $f(t)$ can be defined as the arrival rate of an *effective* boarding demand at time $t$ at the critical stop, where the load of vehicles reaches the maximum. This critical stop will be the beginning stop for one-to-one and one-to-many demand patterns, and the second-to-last stop under many-to-one pattern. The $f$ is a continuous time-varying function and can be obtained by $f(t) \equiv \sum_j f_j(t + \tau_j), t \in [0, E]$, where $f_j$ is the arrival rate of patrons at stop $j$; $\tau_j$ is the travel time of vehicles between the critical stop and stop $j$, which is assumed to be deterministic, fixed, and identical to all buses (imagine the case of bus lines with exclusive lanes); and $E$ denotes the length of the operation period per day.

To serve the above known demand, buses of a fleet with total seats of $M$ are dispatched at a rate of $g(t), t \in [0, E]$, expressed in the unit of seats per hour. For each dispatch, $\lambda$ cost per seat per hour (hours/seat-hour) is related to the operation of a round trip. For each seat in the fleet $M$, a fixed cost $\gamma$ (hours/seat) is associated to the purchase and maintenance. Given the round-trip/cycle time of the bus line being a constant $T$ (hours/cycle), the scheduling problem of shuttle/feeder bus lines is formulated to minimize the total costs for bus operation and patrons' waiting time, as given by

$$\text{minimize}_{g(t),M} \mathcal{J} = \gamma M + \int_0^E [\lambda g(t) + w(t)]dt \tag{1a}$$

subject to:

$$\int_{t-T}^{t} g(s)ds \leq M, t \in [0, E], \tag{1b}$$

$$w(t) = \begin{cases} \frac{cf(t)}{2g(t)}, & \text{if no queue exists at } t \\ \frac{c}{2} + \int_{t_q}^{t} f(s) - g(s)ds, & \text{if queue first forms at } t_q \end{cases}, t \in [0, E], \quad (1c)$$

$$g(t) \geq 0, t \in [0, E], \quad (1d)$$

where (1b) is the fleet size constraint; $w(t)$ (patron-hours/hour) in (1c) is the accumulation rate of patrons' waiting time at $t$; $c$ (seats/bus) is the size/capacity of each bus; $t_q$ is the time when a queue of patrons first forms at the critical bus stop. By 'queue', we mean the phenomenon that the patrons fail to board at least one arrived bus. The last constraint (1d) regulates the valid values of $g(t), t \in [0, E]$. Derivations of (1c) can be found in Hurdle (1973) and omitted here for brevity.

In Hurdle (1973), problem (1) was attacked for peak and off-peak periods, respectively. The peak period was defined as the interval of $[t_q - T, t_d]$, where $t_d$ is the time when the queue dissipates. Other times in $[0, E]$ pertained to the off-period period.

For the peak period, the optimal dispatching rate, $g_{pk}^*(t), t \in [t_q - T, t_d]$, was obtained via graphic analysis, as given by

$$g_{pk}^*(t) = \begin{cases} f(t), & t \in [t_q - T, t_q) \\ g_{pk}^*(t - T), & t \in [t_q, t_d] \end{cases}. \quad (2a)$$

For the off-peak period with no queue, constraint (1b) was ignored by the author, and $w(t)$ in (1c) was reduced to $\frac{cf(t)}{2g(t)}$. Plugging back to (1a), the author then solved the problem via calculus of variations. The optimal off-peak dispatching rate, $g_{op}^*(t)$, was

$$g_{op}^*(t) = \max\left(\sqrt{\frac{cf(t)}{2\lambda}}, f(t)\right), t \in [0, E] \backslash [t_q - T, t_d]. \quad (2b)$$

We note in (2a) and (2b) that the optimized dispatching profiles remain undetermined with $t_q$ and $t_d$ being unknown. To fully resolve the problem, one needs to find the optimal fleet size $M^*$ that determines $t_q$ and $t_d$. The optimal fleet size $M^*$ was, however, not provided in Hurdle (1973). Instead, the author shown the lower and upper bounds of $M^*$, $M_L \leq M^* \leq M_U$, where $M_L$ is defined by that resulting in $t_d - t_q = \gamma$ under the dispatching policy of (2a); and $M_U$ is obtained by the off-peak solution (2b), $M_U \equiv \max_{t \in [0,E]} \int_{t-T}^{t} g_{op}^*(s)ds$. Although trail-and-error methods were suggested to find the optimal $M^*$, the solution efficiency would be low with graphic analysis. This low efficiency may become a serious problem for operators who often have tasks for scheduling hundreds of lines in a city.

In a nutshell, it is found that the scheduling problem (1) still lacks an efficient mathematical solution method to fully resolve the problem.

## 3. Mathematical solutions

This section advances Hurdle's work by presenting a mathematical solution approach to analytically derive the optimal peak-period dispatching rate, $g_{pk}^*(t)$ (Section 3.1). The proposed approach allows us to numerically solve the optimal fleet size $M^*$ (Section 3.2). Notation used in this paper are mostly borrowed from Hurdle (1973) for the sake of consistency. Although phrased in terms of 'patrons', 'buses', and 'bus line/stops', our models apply to other public transportation modes, e.g., subways, trains, airlines, ferries, and ships.

### 3.1. The optimal peak-period dispatching rate

For the pre-queueing interval $t \in [t_q - T, t_q)$, construct a dispatching rate $g(t) = f(t) + k(t)$, where $k(t) \geq 0$ represents an arbitrary path from time $t_q - T$ to $t_q$, and can be interpreted as the departing rate of oversupplied/vacant seats at $t$.[1] Let $\mathcal{T}_q \equiv t_d - t_q$ be the duration that the queue lasts, and $I \equiv \left\lfloor \frac{\mathcal{T}_q}{T} \right\rfloor$ be the number of the cycle time accounted in $\mathcal{T}_q$, where $\lfloor \cdot \rfloor$ returns the closest integer no larger than the argument. The peak-period optimization problem is expressed as

$$\text{minimize}_{g(t),M} \; \mathcal{J}_{pk} = \gamma M + \int_{t_q - T}^{t_d} \lambda g(t) dt + \int_{t_q - T}^{t_d} w(t) dt. \tag{3}$$

subject to: (1b-1d).

Next, we express $\mathcal{J}_{pk}$ in (3) using $k(t)$ and its cumulative function $K(t) \equiv \int_{t_q-T}^{t} k(s) ds$, which means the cumulative number of vacant seats up to $t$. First, the 2nd operation cost-related term at the right-hand-side (RHS) of (3) can be rewritten as

$$\int_{t_q-T}^{t_d} \lambda g(t) dt = \int_{t_q-T}^{t_q} \lambda[f(t) + k(t)] dt + \int_{t_q}^{t_q+\mathcal{T}_q} \lambda g(t) dt. \tag{4}$$

Note that during the queueing interval $[t_q, t_q + \mathcal{T}_q]$, a queue exists and buses are dispatched with zero storage, i.e., the cumulative number of dispatched seats $G(t) \equiv \int_{t_q-T}^{t} g(s) ds$ satisfies the equation of $G(t) = G(t - T) + M$ for $t \in [t_q, t_q + \mathcal{T}_q]$.[2] Thus, we know $g(t) = g(t - T)$ for $t \in [t_q, t_q + \mathcal{T}_q]$. Substituting this result to the third RHS term of (4) yields,

$$\int_{t_q-T}^{t_d} \lambda g(t)dt = \int_{t_q-T}^{t_q} \lambda[f(t) + k(t)]dt + \int_{t_q}^{t_q+\mathcal{T}_q} \lambda g(t - T)dt. \tag{5}$$

---

[1] Hurdle (1973) defined a variable equivalent to the integral of $k(t)$, but unfortunately did not exploit it to analytically solve the problem.

[2] To prove it, consider $\Delta M \geq 0$ seats were hold for an arbitrary interval $[\tilde{t}, \tilde{t} + \Delta t]$ during the queueing period (such that $\tilde{G}(t) = \tilde{G}(t-T) + M - \Delta M \leq \tilde{G}(t-T) + M, t \in [\tilde{t}, \tilde{t} + \Delta t]$) and released later until $\tilde{t} + \Delta \tilde{t}$; then, there would be $\Delta M$ queueing passengers left un-boarded until $\Delta \tilde{t}$ time later. Their waiting time would be increased by $\Delta M \Delta \tilde{t} \geq 0$. The operation cost, however, wasn't changed since $\Delta M$ seats were dispatched after all. Overall, the total system cost would be increased by $\Delta M \Delta \tilde{t} \geq 0$. Thus, $\Delta M$ has to be zero and $G(t) = G(t - T) + M$ must hold for any time during the queueing period, $t \in [t_q, t_q + \mathcal{T}_q]$.

Let $s = t - T$, we can expand the 2$^{nd}$ RHS term of (5) as $\int_{t_q}^{t_q+\mathcal{T}_q} \lambda g(t-T)dt = \int_{t_q-T}^{t_q+\mathcal{T}_q-T} \lambda g(s)ds = \int_{t_q-T}^{t_q} \lambda[f(t) + k(t)]dt + \int_{t_q}^{t_q+\mathcal{T}_q-T} \lambda g(s - T)ds = \cdots$. Thus, (5) can be rewritten as

$$\int_{t_q-T}^{t_d} \lambda g(t)dt = (I + 1)\int_{t_q-T}^{t_q} \lambda[f(t) + k(t)]dt + \int_{t_q-T}^{t_q+\mathcal{T}_q-(I+1)T} \lambda[f(t) + k(t)]ds. \tag{6}$$

where $t_q + \mathcal{T}_q - (I + 1)T \geq t_q - T$ holds since $I \leq \frac{\mathcal{T}_q}{T}$ by definition. Note if $I = \frac{\mathcal{T}_q}{T}$, (6) is reduced to $\int_{t_q-T}^{t_d} \lambda g(t)dt = (I + 1)\int_{t_q-T}^{t_q} \lambda[f(t) + k(t)]dt$, implying the operation cost during the whole peak period is $(I + 1)$ times of that of the pre-queueing interval. This is understandable because the dispatching profile during the queueing period is periodically replicating that of the pre-queueing interval.

In a similar manner, the 3$^{rd}$ waiting-time term in (3) can be expressed in terms of $k(t)$ and $K(t)$ as follows:

$$\int_{t_q-T}^{t_d} w(t)dt = \frac{c}{2}\int_{t_q-T}^{t_q} \frac{f(t)}{f(t)+k(t)}dt + \int_{t_q}^{t_q+\mathcal{T}_q} \left[\frac{c}{2} + F(t) - \left(G(t) - K(t_q)\right)\right]dt, \tag{7}$$

where $F(t) \equiv \int_{t_q-T}^{t} f(s)ds$ is the cumulative number of patrons arrived up to $t$; $K(t_q)$ produces the cumulative vacant seats departed up to $t_q$ (note that no empty seat will depart from the critical stop after $t_q$). Thus, $G(t) - K(t_q)$ returns the cumulative number actual boarded/departed patrons, and $F(t) - \left(G(t) - K(t_q)\right)$ yields the cumulative number of queueing patrons at $t \geq t_q$.

Replacing $G(t)$ by $G(t - T) + M$ in the 2$^{nd}$ RHS term of (7) yields,

$\int_{t_q}^{t_q+\mathcal{T}_q} \left[\frac{c}{2} + F(t) + K(t_q) - G(t - T) - M\right]dt$

$= \int_{t_q-T}^{t_q} \left[\frac{c}{2} + F(t + T) + K(t_q) - G(t) - M\right]dt + \int_{t_q}^{t_q+\mathcal{T}_q-T} \left[\frac{c}{2} + F(t + T) + K(t_q) - G(t - T) - 2M\right]dt$

$= \sum_{i=1}^{2} \int_{t_q-T}^{t_q} \left[\frac{c}{2} + F(t + iT) + K(t_q) - G(t) - iM\right]dt + \int_{t_q}^{t_q+\mathcal{T}_q-2T} \left[\frac{c}{2} + F(t + 2T) + K(t_q) - G(t - T) - 3M\right]dt$

$= \cdots$

$= \sum_{i=1}^{I} \int_{t_q-T}^{t_q} \left[\frac{c}{2} + F(t + iT) + K(t_q) - G(t) - iM\right]dt + \int_{t_q-T}^{t_q+\mathcal{T}_q-(I+1)T} \left[\frac{c}{2} + F(t + (I + 1)T) + K(t_q) - G(t) - (I + 1)M\right]dt.$

By definition $G(t) = F(t) + K(t)$ for $t \in [t_q - T, t_q)$, (7) can then be rewritten as

$$\int_{t_q-T}^{t_d} w(t)dt$$
$$= \frac{c}{2}\int_{t_q-T}^{t_q} \frac{f(t)}{f(t)+k(t)} dt + \sum_{i=1}^{I} \int_{t_q-T}^{t_q} \left[\frac{c}{2} + F(t+iT) + K(t_q) - F(t) - K(t) - iM\right] dt + \quad (8)$$
$$\int_{t_q-T}^{t_q+\mathcal{T}_q-(I+1)T} \left[\frac{c}{2} + F(t+(I+1)T) + K(t_q) - F(t) - K(t) - (I+1)M\right] dt.$$

Using the information from (6) and (8), we define $\mathcal{F}_1(t)$ and $\mathcal{F}_2(t)$, $t \in [t_q - T, t_q)$ as

$$\mathcal{F}_1(t) \equiv \sum_{i=1}^{I+1} \lambda[f(t) + k(t)] + \frac{c}{2}\frac{f(t)}{f(t)+k(t)} + \sum_{i=1}^{I} \left[\frac{c}{2} + F(t+iT) + K(t_q) - F(t) - K(t) - iM\right],$$

$$\mathcal{F}_2(t) \equiv \lambda[f(t) + k(t)] + \left[\frac{c}{2} + F(t+(I+1)T) + K(t_q) - F(t) - K(t) - (I+1)M\right].$$

The original constrained optimization problem (3) can now be rewritten as the following one with respect to $k(t), K(t)$, and $M$ in the pre-queueing interval $[t_q - T, t_q)$:

$$\text{minimize}_{k(t),K(t),M}\; \mathcal{J}_{pk} = \gamma M + \int_{t_q-T}^{t_q+\mathcal{T}_q-(I+1)T} [\mathcal{F}_1(t) + \mathcal{F}_2(t)]dt + \int_{t_q+\mathcal{T}_q-(I+1)T}^{t_q} \mathcal{F}_1(t)dt. \quad (9a)$$

subject to

$$\int_{t_q-T}^{t_q} f(t) + k(t)dt = M, \text{ or } F(t_q) + K(t_q) - F(t_q-T) - K(t_q-T) = M. \quad (9b)$$

And (9) can be further reduced to an equivalent unconstrained variational problem

$$\text{minimize}_{k(t),K(t)}\; \mathcal{J}_{pk} = \int_{t_q-T}^{t_q+\mathcal{T}_q-(I+1)T} j_1(t,k(t),K(t))dt + \int_{t_q+\mathcal{T}_q-(I+1)T}^{t_q} j_2(t,k(t),K(t))dt. \quad (10)$$

where the integrands are

$$j_1(t,k(t),K(t)) = \gamma[f(t)+k(t)] + \sum_{i=1}^{I+2} \lambda[f(t)+k(t)] + \frac{c}{2}\frac{f(t)}{f(t)+k(t)} + \sum_{i=1}^{I+1}\Big[\frac{c}{2} + F(t+iT) +$$
$$K(t_q) - F(t) - K(t) - i[F(t_q) + K(t_q) - F(t_q-T) - K(t_q-T)]\Big], \text{ and}$$

$$j_2(t,k(t),K(t)) = \gamma[f(t)+k(t)] + \sum_{i=1}^{I+1} \lambda[f(t)+k(t)] + \frac{c}{2}\frac{f(t)}{f(t)+k(t)} + \sum_{i=1}^{I}\Big[\frac{c}{2} + F(t+iT) +$$
$$K(t_q) - F(t) - K(t) - i[F(t_q) + K(t_q) - F(t_q-T) - K(t_q-T)]\Big].$$

Problem (10) can be solved using calculus of variations. The Euler-Lagrange equation is

$$\frac{\partial j_1(t,k(t),K(t))}{\partial K(t)} - \frac{d}{dt}\left(\frac{\partial j_1(t,k(t),K(t))}{\partial k(t)}\right) = 0 \text{ for } t \in [t_q-T, t_q+\mathcal{T}_q-(I+1)T], \text{ and}$$

$$\frac{\partial j_2(t,k(t),K(t))}{\partial K(t)} - \frac{d}{dt}\left(\frac{\partial j_2(t,k(t),K(t))}{\partial k(t)}\right) = 0 \text{ for } t \in t \in (t_q+\mathcal{T}_q-(I+1)T, t_q),$$

which become

$$0 - \frac{d}{dt}\left(\gamma + \sum_{i=1}^{I+2} \lambda + \frac{c}{2}\frac{\partial}{\partial k(t)}\left(\frac{f(t)}{f(t)+k(t)}\right)\right) = 0 \text{ for } t \in [t_q-T, t_q+\mathcal{T}_q-(I+1)T], \quad (11a)$$

$$0 - \frac{d}{dt}\left(\gamma + \sum_{i=1}^{I+1} \lambda + \frac{c}{2}\frac{\partial}{\partial k(t)}\left(\frac{f(t)}{f(t)+k(t)}\right)\right) = 0 \text{ for } t \in (t_q+\mathcal{T}_q-(I+1)T, t_q). \quad (11b)$$

The solution to (11a, b) is

$$g^*_{pk}(t) = \begin{cases} \max\left(\sqrt{\frac{c}{2}\frac{f(t)}{\gamma+(I+2)\lambda}}, f(t)\right), & t \in [t_q - T, t_q + \mathcal{T}_q - (I+1)T] \\ \max\left(\sqrt{\frac{c}{2}\frac{f(t)}{\gamma+(I+1)\lambda}}, f(t)\right), & t \in (t_q + \mathcal{T}_q - (I+1)T, t_q) \end{cases}. \quad (12)$$

The optimal peak-period dispatching rate in (12) has three special forms. The first is Hurdle's solution (2a) when the peak demand is high (with $\gamma$ and/or $\lambda$ being not too small) such that $g^*_{pk}(t) = f(t)$ holds for $t \in [t_q - T, t_q)$; the second is in the form of the off-peak solution $g^*_{pk}(t) = \max\left(\sqrt{\frac{c}{2}\frac{f(t)}{\lambda}}, f(t)\right)$, $t \in [0, E]$, when $\gamma = 0$ and there is no fleet size constraint (implying no queue and $I = 0$); and the third is $g^*_{pk}(t) = \max\left(\sqrt{\frac{c}{2}\frac{f(t)}{\gamma}}, f(t)\right)$, $t \in [t_q - T, t_q)$ when $\lambda = 0$.

Among the above three, Hurdle's result is of the most significance from a practical perspective. This is because realistic conditions oftentimes suffice $f(t) \geq \sqrt{\frac{c}{2}\frac{f(t)}{\gamma+(I+2)\lambda}} \geq \sqrt{\frac{c}{2}\frac{f(t)}{\gamma+(I+1)\lambda}} \Rightarrow f(t) \geq \frac{c}{2}\frac{1}{\gamma+(I+1)\lambda} > \frac{c}{2}\frac{1}{\gamma+(I+2)\lambda}$ during peak period. To see this fact, consider ordinary busy bus lines where the arrival rate of patrons can easily surpass dozens per mins; and conservatively, suppose an arrival rate of $f(t_q - T) = 20$ (patrons/min) and a bus size of $c = 40$ (patrons/bus). It can be seen that $20 \geq \frac{40}{2}\frac{1}{\gamma+(I+1)\lambda} > \frac{40}{2}\frac{1}{\gamma+(I+2)\lambda}$ holds for realistic values of $\gamma$ and $\lambda$, which are much larger than 1 (mins/seat or mins/seat-hour) (Sivakumaran et al., 2014; Gu et al., 2016; Chang and Schonfeld, 1991; Chen et al., 2020). Similar arguments also apply to other transit modes, e.g., subways and Bus Rapid Transit (BRT). Interested readers can verify the above inequality relationship with the corresponding mode-dependent parameter values.

### 3.2. The optimal fleet size

The formulation of (9a) allows us to find the first-order condition with respect to the optimal fleet size $M^*$,

$$\frac{d\mathcal{J}_{pk}}{dM} = \gamma - \int_{t_q-T}^{t_q} \sum_{i=1}^{I} i\, dt - \int_{t_q-T}^{t_q+\mathcal{T}_q-(I+1)T}(I+1)dt = 0,$$

which can be reorganized as

$$\gamma = \frac{I(I+1)}{2}T + (I+1)(\mathcal{T}_q - IT), \quad (13)$$

where $I = \left\lfloor \frac{\mathcal{T}_q}{T} \right\rfloor$ by definition. If $I = \frac{\mathcal{T}_q}{T}$, (13) will be reduced to $\gamma = \frac{I(I+1)}{2}T$.

The solution of $\mathcal{T}_q$ (and $I$) to (13) can be used to find the $M^*$. To do so, we use the following two boundary conditions at time $t_q$ and $t_q + \mathcal{T}_q$, respectively:

$$G(t_q) - G(t_q - T) = M, \quad (14)$$

$$F(t_q + \mathcal{T}_q) = G(t_q + \mathcal{T}_q). \quad (15)$$

In (15), $G(t_q + \mathcal{T}_q)$ can be folded into the pre-queueing interval by $G(t_q + \mathcal{T}_q) = G(t_q + \mathcal{T}_q - T) + M = G(t_q + \mathcal{T}_q - 2T) + 2M = \cdots = G(t_q + \mathcal{T}_q - IT) + IM$, which can be further reformulated as

$$F(t_q + \mathcal{T}_q) = \int_{t_q}^{t_q+\mathcal{T}_q-IT} g(t)dt + (I+1)M. \quad (16)$$

Substituting (14) in (16) and replacing $g(t)$ by $g(t-T)$ at interval $[t_q, t_q + \mathcal{T}_q - IT]$ yield,

$$F(t_q + \mathcal{T}_q) = \int_{t_q-T}^{t_q+\mathcal{T}_q-(I+1)T} g(t)dt + (I+1)\left(G(t_q) - G(t_q - T)\right). \quad (17)$$

Given $\mathcal{T}_q$ from (13) and $g$ and $G$ taking the optimal solution $g_{pk}^*$ and $G_{pk}^*$, we now have two equations (14, 17) and two unknown variables $t_q, M$, which can be solved. The solution may be remarkably straightforward. Typically, we have $\gamma \leq T$ in reality; then (13) is reduced to $\mathcal{T}_q = \gamma$ with $I = 0$;[3] and (17) is simplified to $F(t_q + \gamma) = G(t_q + \gamma - T) + G(t_q) - 2G(t_q - T)$, of which the solution of $t_q^*$ is plugged into (14) yielding the optimal fleet size $M^*$. For other cases of $\gamma > T$, $\mathcal{T}_q^*, I^*|\gamma$ should be first solved by (13), and then $t_q^*|\mathcal{T}_q^*, I^*$ by (17), and lastly, the optimal fleet size $M^*|t_q^*$.

Ultimately, the fleet size in terms of actual vehicles can be determined by $\left[\frac{M^*}{c}\right]$, where $[\cdot]$ rounds the argument to the closest integer.

With above solutions, detailed system performance metrics can then be obtained, including the beginning and ending times of queueing $(t_q, t_d)$, the total waiting time in queues ($\int_{t_q}^{t_d} w(t)dt$), the maximum queue length ($\max_{t\in[t_q,t_d]}[F(t) - G(t)]$), etc. These metrics would otherwise be difficult to gain via graphic analysis and trial-and-error method.

### 3.3. The optimal off-peak dispatching rate

As shown in Section 3.1, the cost parameters, $\gamma$ and $\lambda$, have theoretical influences on the optimal solution to the peak-period dispatching problem. This is also true for the off-peak problem. For instance, when $\lambda = 0$ the off-peak dispatching rate, as predicted by Hurdle's formula (2b), tends to be infinitely large. This is obviously impossible given the fleet size constraint. Therefore, a general off-peak dispatching problem should be written as

$$\text{minimize}_{g(t)}\ \mathcal{J}_{op} = \int_0^E \left[\lambda g(t) + \frac{cf(t)}{2g(t)}\right] dt \tag{18a}$$

subject to:

$$\int_{t-T}^{t} g(s)ds \leq M, t \in [0, E]\setminus[t_q - T, t_d], \tag{18b}$$

$$g(t) \geq f(t), t \in [0, E]\setminus[t_q - T, t_d]. \tag{18c}$$

where the fleet size $M$ is determined by the peak-period dispatching problem (section 3.2). The above problem (18) is a variational problem with inequality constraint.[4] In general, the solution would be two dispatching policies that alternate during the off-peak period: (i) $g_{op}^*(t) = \max\left(\sqrt{\frac{cf(t)}{2\lambda}}, f(t)\right)$ for some subintervals in $[0, E]\setminus[t_q - T, t_d]$ when (18b) is inactive; and otherwise (ii) $g_{op}^*(t) = g_{op}^*(t - T)$. To determine the beginning and ending times of these subintervals, detailed solution procedures can be found in Salzborn (1972) and Kamien and Schwartz (1991), and thus omitted here for the sake of brevity.

---

[3] To see why, consider the fact that $\frac{I(I+1)}{2}T \geq T$ for $I = 1,2,\ldots$ and $\mathcal{T}_q - IT \geq 0$ by definition; thus, to satisfy $\gamma = \frac{I(I+1)}{2}T + (I+1)(\mathcal{T}_q - IT) \leq T, I = 0$ must hold. Therefore, (12) is reduced to $\mathcal{T}_q = \gamma$ for $\gamma \leq T$.

Fortunately, the above solutions will be reduced to Hurdle's result (2b) if constraint (18b) is always satisfied during off-peak period. This special case will happen more than often in reality when $M$ is sufficiently large (in other words, the peak demand is high) and $\lambda$ is not too small. To see this, consider again bus lines with $c = 40$ (seats/bus) and $\lambda \gg 1$ mins/seat-hour and conservatively suppose the peak period begins as early as when $f(t)$ reaches 20 patrons/min. Obviously, the inequality of $\max_{t \in [0,E]\setminus[t_q-T,t_d]} \int_{t-T}^{t} g_{op}^*(s)ds < 20T < \int_{t_q-T}^{t_q} f(s)ds = M$ holds true and remains so for busier lines.

### 3.4. Numerical examples

For the sake of illustration, we make up the following time-varying arrival pattern of patrons:

$$f(t) = D \cdot Tr\mathcal{N}(\mu, \sigma^2, 0, E), \tag{19}$$

where $D$ is the total trips during the study period of $E$ hours; $\mu, \sigma$ denote the mean and standard deviation of the normal distribution of the number of arrivals; and $Tr\mathcal{N}(\cdot,\cdot,0,E)$ represents the probability density function of the normal distribution truncated by the interval $[0, E]$. Demand of (19) has one peak (at time $\mu$) that is used to represents a morning or evening peak.

---

[4] A similar constrained scheduling problem was studied by Salzborn (1972).

Select values of the parameters used in the above models are given in Table 2. Two values of the fleet-cost parameter $\gamma$ (i.e., 30 and 90 mins/seat as opposed to $T = 60$ mins) are chosen to denote a normal and a costly-fleet scenario, respectively. Other parameter values are set to generally represent the realistic situations.

**Table 2. Parameter values.**

| Parameters | Baseline values | Parameters | Baseline values |
|---|---|---|---|
| $D$ (patrons) | 1000 | $E$ (mins) | 180 |
| $T$ (mins) | 60 | $c$ (patrons/bus) | 25 |
| $\mu$ (mins) | 60 | $\sigma$ (mins) | 30 |
| $\gamma$ (mins/seat) | {30, 90} | $\lambda$ (mins/seat-hour) | 5 |

With the above baseline settings, the optimized dispatching rates coincide with Hurdle's results, i.e., $g_{pk}^*(t) = f(t), t \in [t_q - T, t_q)$ and $g_{op}^*(t) = \max\left(\sqrt{\frac{cf(t)}{2\lambda}}, f(t)\right), t \in [0, E]\setminus[t_q - T, t_q)$. Figures 1a and 1b depicts the dispatching profiles and the cumulative number of dispatches, respectively, under the normal scenario with $\gamma = 30$ mins/seat. As seen, the queueing time starts at 75 mins and lasts $\mathcal{T}_q = 30$ mins. During the time, the optimized fleet size constraint with $M^* = 640$ (seats) binds. Figures 1c and 1d are for the costly-fleet scenario. It is observed that the queueing starting time moves forward to 59 mins and the duration extends to $\mathcal{T}_q = 75$ mins, which covers more than one $T$ periods. In this case, the optimized feet size becomes $M^* = 477$ seats.

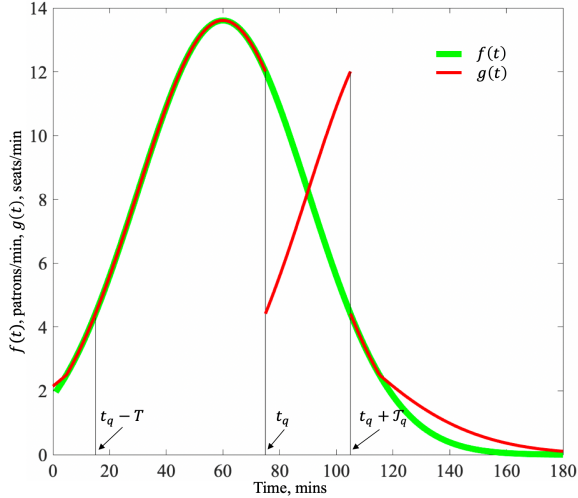

a. Dispatching rate with $\gamma = 30$ mins/seat

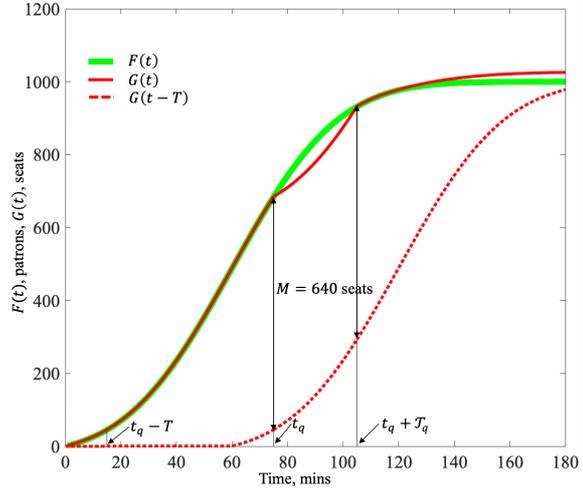

b. Cumulative dispatches with $\gamma = 30$ mins/seat

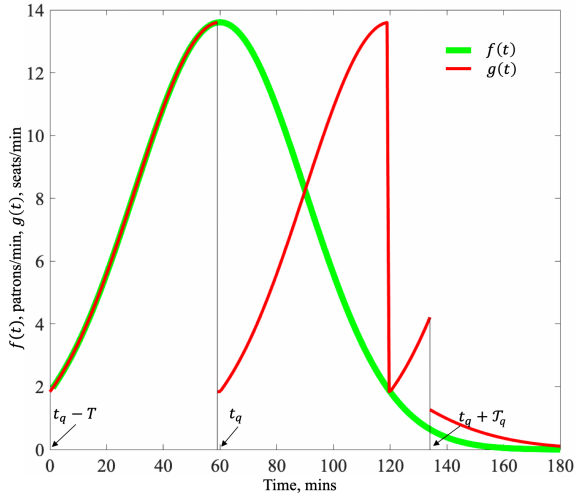 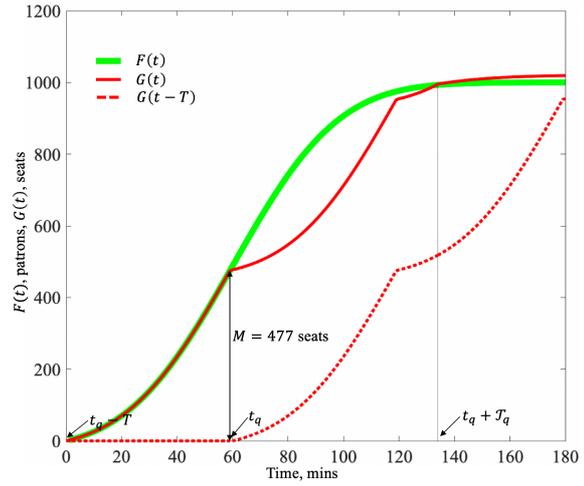

c. Dispatching rate with $\gamma = 90$ mins/seat

d. Cumulative dispatches with $\gamma = 90$ mins/seat

**Figure 1. Visualization of the optimal dispatching rates.**

The system performance can be measured by several metrics, as exemplified in Table 3. For instance, the total waiting times in queue under the above two scenarios are 1177.8 mins and 9131.8 mins, respectively; and the maximum queue length are 59 patrons and 223 patrons, respectively. These metrics would be useful for operators to monitor the system and accordingly develop necessary countermeasures, e.g., queue management.

**Table 3. System performance metrics.**

| Metrics<br>Scenarios | Waiting time in queue (mins) | | Number of queueing patrons | |
|---|---|---|---|---|
| | Total | Average | Maximum | Average |
| Scenario I ($\gamma = 30$) | 1177.8 | 4.76 | 59 | 20 |
| Scenario II ($\gamma = 90$) | 9131.8 | 17.67 | 223 | 152 |

Thanks to the obtained optimal conditions, the solution finding is extremely computationally efficient. Such a high efficiency allows us to conduct parameter analysis, which might be unpractical by graphic analysis. Figure 2 shows the results of the changes in $M^*$ and $t_q$ with respect to $\gamma \in [5,50]$ mins/seat. While both decrease with the rising $\gamma$ values, the queueing starting time $t_q$ appears to diminish in a linear manner. These findings are not available in the literature, to the best of our knowledge.

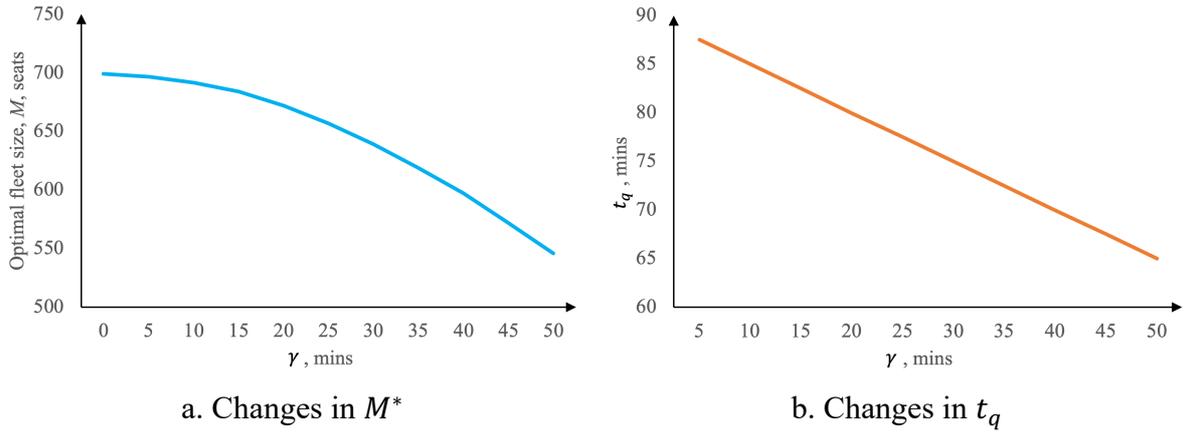

a. Changes in $M^*$

b. Changes in $t_q$

**Figure 2. Changes of the optimal solutions with respect to varying $\gamma$.**

Sensitivity analysis is also done to $\sigma$ that represents the temporal aggregation of patrons' arrivals. As expected, smaller/larger $\sigma$ values, indicating high/low-level aggregation, are found to require more and fewer fleet sizes, respectively.

We also examined scenarios with very low demand levels and confirmed that our optimal peak-period dispatching rate diverges from and outperforms Hurdle's. For instance, when $D = 50$ (patrons) and $f(t) < 0.7$ (patrons/min) for $t \in [0, E]$, the two dispatching profiles are shown in Figures 3a and 3b. As seen in Figure 3b, $g_{pk}^*(t) > f(t)$ holds for a short interval of the peak period. Not surprisingly, the optimized system cost is less than that under Hurdle's dispatching policy; see the 2nd column of Table 4. Table 4 also summarizes other system metrics including the Total Waiting Time in Queue (TWTQ), the queueing beginning time, and the fleet size. As seen, the proposed optimal dispatching policy also leads to less TWTQ. As a result of the increased dispatching rate, the queueing beginning time shifts earlier than that of Hurdle's policy. The fleet size is not affected by such a small divergence.

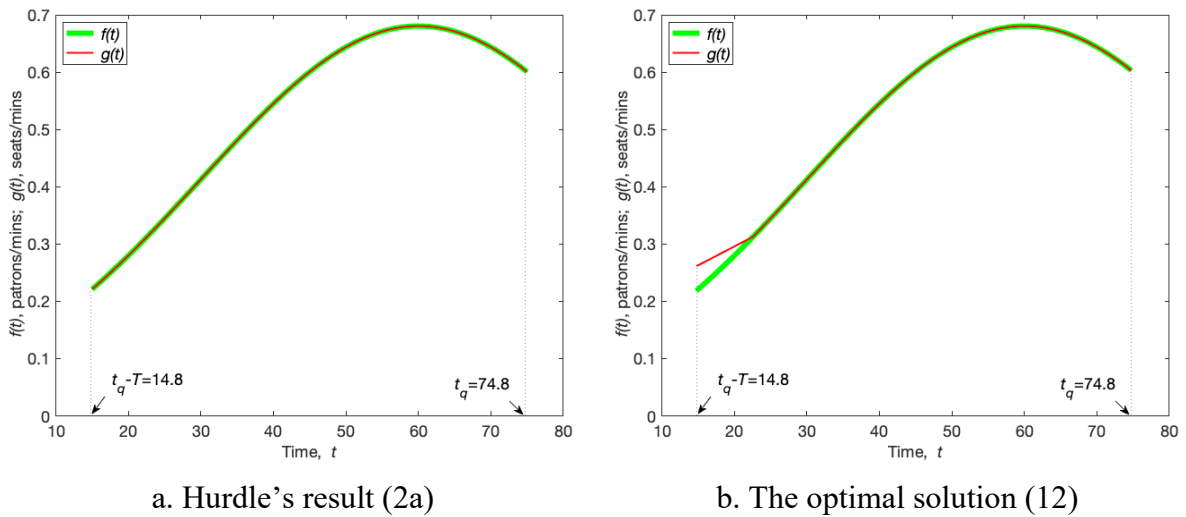

a. Hurdle's result (2a)

b. The optimal solution (12)

**Figure 3. Dispatching profiles of Hurdle's and ours under low demand level.**

Table 4. Results comparison under low demand level.

| Dispatching rates | System cost, $J_{pk}$ (mins) | TWTQ (mins) | Queueing beginning time, $t_q$ | Fleet size, $M$ (seats) |
|---|---|---|---|---|
| Hurdle's | 2364 | 58.9 | 75 | 32 |
| Ours | 2358 | 56.7 | 74.8 | 32 |

The above results demonstrate that Hurdle's results suffice to apply to moderately and heavily busy transit lines. For light transit lines, the arrival rate of patrons might be too low for CA models to apply. If so, mathematical models, e.g., those in Newell (1971) and Sheffi and Sugiyama (1982), should be used instead.

## 4. Two extensions

This section demonstrates two possible extensions to Hurdle's scheduling problem. The first is about general transit lines with multiple origins and destinations (Section 4.1); and the second considers mixed-size/modular vehicles (Section 4.2). Without particular notification, we will be concerned about busy transit lines where fleet size would be sufficiently large and the corresponding constraint won't bind during the off-peak period.

### 4.1. Scheduling for transit lines with multiple origins and destinations

Consider a bus line with multiple stops indexed by $j \in \mathbf{J} \equiv \{0,1,...,J\}$, where $J$ is the total number of stops. Patrons originate and destinate between any pairs of stops. Let the arrival rate of patrons at stop $j$ at time $t + \tau_j$ be $f_j(t) \equiv f_j(t + \tau_j)$, where $\tau_j$ is the average travel time of buses traversing from the depot (indexed by 0 and $J$) to the $j^{th}$ stop. Define $\beta_{od}(t) \equiv \frac{f_{od}(t)}{f_o(t)}$ (patrons/hour) as the ratio of demand traveling from stop $o$ to stop $d$, where $f_{od}(t)$ is the demand between the stop pair, $o,d \in \mathbf{J}$. Then, the waiting time at any stop $j \in \mathbf{J}$ can be formulated as

$$w_j(t) = \begin{cases} \frac{cf_j(t)}{2g(t)}, & \text{if no queue exists at stop } j \text{ for disptach time } t \\ \frac{c}{2} + \int_{t_{j,q}}^{t}[f_j(s) - b_j(s)]ds, & \text{if queue first forms at } t_{j,q} \end{cases}, j \in \mathbf{J}, t \in [0,E], \quad (20)$$

where $t_{j,q}$ is the time that a queue first forms at stop $j$; $b_j(s)$ (patrons/hour) is the actual boarding rate at stop $j$ for the dispatch at time $s$, which is given by

$$b_j(s) = \min\{e_j(s), f_j(s)\}, \tag{21a}$$

$$e_j(s) = \frac{g(s)}{c}\left[c - \frac{c}{g(s)}\left(\sum_{o=0}^{j-1} b_o(s) - \sum_{o=0}^{j-1}\sum_{d=o+1}^{j} b_o(s)\beta_{od}(s)\right)\right], \tag{21b}$$

where $e_j(s)$ in (21a) denotes the arrival rate of available seats at stop $j$, and is defined by (21b); min{·} in (21a) returns the minimum value of the arguments, i.e., the actual boarding rate per hour. In (21b), $\sum_{o=0}^{j-1} b_o(s)$ and $\sum_{o=0}^{j-1} \sum_{d=o+1}^{j} b_o(t)\beta_{od}(t)$ (patrons/hour) mean the cumulative boarding patrons from stop 0 to $j-1$ and the cumulative alighting ones up to stop $j$, respectively; the term in bracket means the available seats per arriving vehicle, multiplying by $\frac{g(t)}{c}$ (i.e., the arrival rate of vehicles) yields the definition of $e_j(s)$.

With (20, 21) in hand, the scheduling problem of general bus lines is rewritten as

$$\text{minimize}_{g(t),M} \; \mathcal{J} = \gamma M + \int_0^E \left[\lambda g(t) + \sum_0^J w_j(t)\right] dt \quad (22)$$

subject to: (1b, 1d), and (20, 21).

Next, the off-peak scheduling problem is first examined (subsection 4.1.1) and then the peak-period problem (subsection 4.1.2), followed by the optimal fleet size problem (subsection 4.1.3).

### 4.1.1. The optimal off-peak dispatching rate

Define the off-peak period $\mathbf{E}_{op}$ as time when the fleet size constraint is not binding, and we have

**Proposition 1**. During the off-peak period, there is either no queue at any stops along the line or only short-lasting queue(s).

**Proof.** The average waiting time in queue per patron, if not zero, must be less than the operation cost per seat dispatch, $\lambda$; otherwise, additional seat(s) will be dispatched immediately from the pool of idling vehicles to prevent the queue from forming. Thus, the queue(s), if exists, must be short-time lasting. ∎

To facilitate the modeling, we assume

**Assumption 1.** The $\lambda$ is small-valued such that any queue will be prevented from forming during the off-peak period but the fleet size constraint is not activated.

This assumption is not strict and very much in line with reality; consider how rare the case when patrons cannot board the first-arrival bus and the bus storage reaches zero during off-peak hours. Per Assumption 1, the off-peak scheduling problem is modeled as

$$\text{minimize}_{g(t)} \; \mathcal{J}_{op} = \int_{t \in \mathbf{E}_{op}} \left[\lambda g(t) + \sum_{j=0}^{J} \frac{c f_j(t)}{2g(t)}\right] dt \quad (23a)$$

subject to:

$$g(t) \geq r_m(t) \equiv \max_{j \in \mathbf{J}}\left(r_j(t)\right), t \in \mathbf{E}_{op}, \quad (23b)$$

where (23b) guarantees that no queue forms at any stops during the off-peak period; $r_m(t)$ denotes the maximum cross-sectional departing flux along the line, and $r_j(t) \equiv \left(\sum_{o=0}^{j} f_o(t) - \sum_{o=0}^{j-1} \sum_{d=o+1}^{j} f_o(t)\beta_{od}(t)\right)$ is the cross-sectional demand flux departing from stop $j$.

According to the first-order condition of (23), i.e., $\frac{d\mathcal{J}_{op}}{dg(t)} = \lambda - \frac{c\sum_{j=0}^{J}f_j(t)}{2g^2(t)} = 0$ for every $t \in \mathbf{E}_{op}$, we obtain the optimal off-peak dispatching rate as (24) that combines the boundary constraint (23b).

$$g_{op}^*(t) = \max\left\{\sqrt{\frac{c}{2\lambda}\sum_{j=0}^{J}f_j(t)}, r_m(t)\right\}, t \in \mathbf{E}_{op}. \tag{24}$$

It is easy to verify that (24) will be reduced to (2b) when the transit line is of shuttle or feeder services, i.e., $r_m(t) = f(t) \equiv \sum_{j=0}^{J}f_j(t)$.

### 4.1.2. The optimal peak-period dispatching rate

Now, let's turn our attention to the peak-period scheduling problem. From a practical point of view, we assume the following dispatching principle.

**Assumption 2**. When queues form at any stops of the bus line, buses will be dispatched as soon as they are available.[5]

Per Assumption 2, define the peak period as $[t_c - T, t_e]$, where $[t_c, t_e]$ is the *critical period* when the fleet size constraint is binding, i.e., $G(t) = G(t - T) + M$ and $g(t) = g(t - T), t \in [t_c, t_e]$. The $t_c$ and $t_e$ are variables to be determined later in association with the optimal fleet size $M^*$.

Next, we first derive the total waiting times at all stops during the peak period, i.e., $\int_{t_c-T}^{t_e}\sum_{0}^{J}w_j(t)\,dt$. In doing so, stops are classified into three groups according to their states: (i) undersaturated stops without queue at $t \in [t_c - T, t_e]$, denoted $u_t \in \mathbf{J}$; (ii) oversaturated stops with partial boarding demand left in queue at $t \in [t_c - T, t_e]$, $p_t \in \mathbf{J}$, and (iii) congested stops with all boarding demand left in queue at $t \in [t_c - T, t_e]$, $q_t \in \mathbf{J}$. At each undersaturated stop, the total waiting time is straightforward, i.e., $w_{u_t}(t) = \frac{c}{2}\frac{f_{u_t}(t)}{g(t)}, u_t \in \mathbf{J}$. For the latter two types of stops, define $r_j^a(t) \equiv r_j(t) - f_j(t)$ meaning the arriving passenger flux before boarding at stop $j \in \{p_t, q_t\}$. Then, $\max(g(t) - r_j^a(t), 0)$ yields the arrival rate of available seats for boarding demand. Thus, the incremental rate of queueing patrons is

$$f_j(t) - \max(g(t) - r_j^a(t), 0) = \begin{cases} r_j(t) - g(t), \text{if } g(t) > r_j^a(t), j = p_t \\ f_j(t), \text{if } g(t) \leq r_j^a(t), j = q_t \end{cases}.[6]$$

---

[5] Admittedly, such a dispatching principle may not be the optimum. Buses could be hold to serve more boarding patrons at some stops by causing queue(s) at certain stop(s) last slightly longer. However, we rule out this possibility for the merits of the following analytical analysis.
[6] This result is exact for the first queueing (oversaturated or congested) stop $j$ with no queue at any upstream stops. For non-first queueing stops, $j'$, the result is conservative as if the through demand at $j'$ had a higher boarding priority at upstream queueing stops than the alighting demand at $j'$.

Therefore, the total waiting time over the peak period $[t_c - T, t_e]$ at all stops can be expressed by the sum of three terms, i.e., the total waiting time at three types of stops, which are given by

$$\int_{t_c-T}^{t_e} \sum_{j=1}^{J} w_j(t)\,dt = \left[\int_{t_c-T}^{t_e} \left(\frac{c}{2}\frac{\Sigma u_t f_{u_t}(t)}{g(t)}\right) dt\right] + \left[\int_{t_c-T}^{t_e} \Sigma p_t \left(\frac{c}{2} + R_{p_t}(t) - \left(G(t) - \mathcal{K}_{p_t}\right)\right) dt\right] + \left[\int_{t_c-T}^{t_e} \Sigma q_t \left(\frac{c}{2} + F_{q_t,q}(t)\right) dt\right], \quad (25)$$

where $R_{p_t}(t) \equiv \int_{t_c-T}^{t} r_{p_t}(s)ds$ is the cumulative function of arriving flux; $\mathcal{K}_{p_t} \equiv G(t_{p_t,q}) - R_{p_t}(t_{p_t,q})$ means the cumulative vacant seats departed up to the time when the queue first forms at stop $p_t$ at time $t_{p_t,q} \in [t_c - T, t_e]$; and $F_{q_t,q}(t) \equiv \int_{t_{q_t,q}}^{t} f_{q_t}(s)ds$ represents the cumulative number of queueing patrons at stop $q_t$ since the time the queue first forms, i.e., $t_{q_t,q}$.

Then, in a similar manner as in Section 3.1, we construct a dispatching profile for the pre-critical period, $g(t) = r_u(t) + k(t), t \in [t_c - T, t_c)$, where $r_u(t) \geq 0$ is an unknown base function and $k(t) \geq 0$ a decision function to be determined. The relationship between the corresponding cumulative functions is $G(t) = R_u(t) + K(t), t \in [t_c - T, t_c)$.

Replacing $g(t)$ by $r_u(t) + k(t)$, $G(t)$ by $R_u(t) + K(t)$ for $t \in [t_c - T, t_c)$, and $G(t)$ by $G(t-T) + M$ and $g(t) = g(t-T)$ for $t \in [t_c, t_e]$, (25) can be reorganized as

$$\int_{t_c-T}^{t_e} \sum_{j=1}^{J} w_j(t)\,dt = \left[\int_{t_c-T}^{t_c} \left(\frac{c}{2}\frac{\Sigma u_t f_{u_t}(t)}{r_u(t)+k(t)}\right) dt + \int_{t_c}^{t_e} \left(\frac{c}{2}\frac{\Sigma u_t f_{u_t}(t)}{g(t-T)}\right) dt\right]$$

$$+ \left[\int_{t_c-T}^{t_c} \Sigma p_t \left(\frac{c}{2} + R_{p_t}(t) - R_u(t) - K(t) + \mathcal{K}_{p_t}\right) dt + \int_{t_c}^{t_e} \Sigma p_t \left(\frac{c}{2} + R_{p_t}(t) - G(t-T) - M + \mathcal{K}_{p_t}\right) dt\right]$$

$$+ \left[\int_{t_c-T}^{t_c} \Sigma q_t \left(\frac{c}{2} + F_{q_t,q}(t)\right) dt + \int_{t_c}^{t_e} \Sigma q_t \left(\frac{c}{2} + F_{q_t,q}(t)\right) dt\right]$$

$$= \left[\int_{t_c-T}^{t_c} \left(\frac{c}{2}\frac{\Sigma u_t f_{u_t}(t)}{r_u(t)+k(t)}\right) dt + \int_{t_c-T}^{t_e-T} \left(\frac{c}{2}\frac{\Sigma u_{(s+T)} f_{u_{(s+T)}}(s+T)}{g(s)}\right) ds\right]$$

$$+ \left[\int_{t_c-T}^{t_c} \Sigma p_t \left(\frac{c}{2} + R_{p_t}(t) - R_u(t) - K(t) + \mathcal{K}_{p_t}\right) dt + \int_{t_c-T}^{t_e-T} \Sigma p_{(s+T)} \left(\frac{c}{2} + R_{p_{(s+T)}}(s+T) - G(s) - M + \mathcal{K}_{p_{(s+T)}}\right) ds\right]$$

$$+ \left[\int_{t_c-T}^{t_c} \Sigma q_t \left(\frac{c}{2} + F_{q_t,q}(t)\right) dt + \int_{t_c-T}^{t_e-T} \Sigma q_{(s+T)} \left(\frac{c}{2} + F_{q_{(s+T)},q}(s+T)\right) ds\right]$$

$$= \left[\int_{t_c-T}^{t_c} \left(\frac{c}{2}\frac{\Sigma u_t f_{u_t}(t)}{r_u(t)+k(t)}\right) dt + \int_{t_c-T}^{t_c} \left(\frac{c}{2}\frac{\Sigma u_{(s+T)} f_{u_{(s+T)}}(s+T)}{r_u(s)+k(s)}\right) ds + \int_{t_c}^{t_e-T} \left(\frac{c}{2}\frac{\Sigma u_{(s+T)} f_{u_{(s+T)}}(s+T)}{g(s-T)}\right) ds\right]$$

$$+ \left[\int_{t_c-T}^{t_c} \Sigma p_t \left(\frac{c}{2} + R_{p_t}(t) - R_u(t) - K(t) + \mathcal{K}_{p_t}\right) dt + \int_{t_c-T}^{t_c} \Sigma p_{(s+T)} \left(\frac{c}{2} + R_{p_{(s+T)}}(s+T) - R_u(s) - K(s) - M + \mathcal{K}_{p_{(s+T)}}\right) ds + \int_{t_c}^{t_e-T} \Sigma p_{(s+T)} \left(\frac{c}{2} + R_m(s+T) - G(s-T) - 2M + \mathcal{K}_{p_{(s+T)}}\right) ds\right]$$

$$+ \left[\int_{t_c-T}^{t_c} \Sigma q_t \left(\frac{c}{2} + F_{q_t,q}(t)\right) dt + \int_{t_c-T}^{t_c} \Sigma q_{(s+T)} \left(\frac{c}{2} + F_{q_{(s+T)},q}(s+T)\right) ds + \int_{t_c}^{t_e-T} \Sigma q_t \left(\frac{c}{2} + F_{q_t,q}(s+T)\right) ds\right]$$

$$= \cdots$$

$$\begin{aligned}
&= \left[ \int_{t_c-T}^{t_c} \left( \frac{c}{2} \frac{\sum_{i=0}^{I} \sum u_{(t+iT)} f_{u_{(t+iT)}}(t+iT)}{r_u(t)+k(t)} \right) dt + \int_{t_c-T}^{t_e-(I+1)T} \left( \frac{c}{2} \frac{\sum u_{(t+IT+T)} f_{u_{(t+IT+T)}}(t+IT+T)}{r_u(t)+k(t)} \right) ds \right] \\
&+ \left[ \int_{t_c-T}^{t_c} \sum_{i=0}^{I} \sum p_{(t+iT)} \left( \frac{c}{2} + R_{p_{(t+iT)}}(t+iT) - R_u(t) - K(t) - iM + \mathcal{K}_{p_{(t+iT)}} \right) dt + \right. \\
&\left. \int_{t_c-T}^{t_e-(I+1)T} \sum p_{(t+IT+T)} \left( \frac{c}{2} + R_{p_{(t+IT+T)}}(t+(I+1)T) - R_u(t) - K(t) - (I+1)M + \mathcal{K}_{p_{(t+IT+T)}} \right) ds \right] \\
&+ \left[ \int_{t_c-T}^{t_c} \sum_{i=0}^{I} \sum q_{(t+iT)} \left( \frac{c}{2} + F_{q_t,q}(t+iT) \right) dt + \int_{t_c-T}^{t_e-(I+1)T} \sum q_{(t+IT+T)} \left( \frac{c}{2} + F_{q_{(t+IT+T)},q}(t+IT+T) \right) ds \right].
\end{aligned} \quad (26)$$

where $\mathcal{K}_{p_{(t+iT)}} = G\left(t_{p_{(t+iT)},q}\right) - R_{p_t}\left(t_{p_{(t+iT)},q}\right)$ can be expressed as $\mathcal{K}_{p_{(t+iT)}} = R_u(t_{p,q}) + K(t_{p,q}) - R_{p_t}(t_{p,q} + i_{p,q}T) + i_{p,q}M$ by letting $t_{p_{(t+iT)},q} = t_{p,q} + i_{p,q}T, t_{p,q} \in [t_c - T, t_c), 0 \le i_{p,q} \le i$.

We note that the peak-period operation cost term, $\int_{t_c-T}^{t_e} \lambda g(t)$, can be similarly treated as in Section 3.1. And consequently, an unconstrained peak-period scheduling problem can be constructed with the Euler-Lagrange equation being

$$\frac{\partial \sum_{j \in J} w_j(t)}{\partial K(t)} - \frac{d}{dt} \left[ \gamma + (I+2)\lambda + \frac{\partial \sum_{j \in J} w_j(t)}{\partial k(t)} \right] = 0, t \in [t_c - T, t_e - IT - T],$$

$$\frac{\partial \sum_{j \in J} w_j(t)}{\partial K(t)} - \frac{d}{dt} \left[ \gamma + (I+1)\lambda + \frac{\partial \sum_{j \in J} w_j(t)}{\partial k(t)} \right] = 0, t \in (t_e - IT - T, t_c) \Rightarrow$$

$$0 - \frac{d}{dt} \left[ \frac{\partial}{\partial k(t)} \left( \gamma + (I+2)\lambda + \frac{c}{2} \frac{\sum_{i=0}^{I+1} \sum u_{(t+iT)} f_{u_{(t+iT)}}(t+iT)}{r_u(t)+k(t)} \right) \right] = 0, t \in [t_c - T, t_e - IT - T], \quad (27a)$$

$$0 - \frac{d}{dt} \left[ \frac{\partial}{\partial k(t)} \left( \gamma + (I+1)\lambda + \frac{c}{2} \frac{\sum_{i=0}^{I} \sum u_{(t+iT)} f_{u_{(t+iT)}}(t+iT)}{r_u(t)+k(t)} \right) \right] = 0, t \in (t_e - IT - T, t_c), \quad (27b)$$

The solution to (27a, b) is

$$\tilde{g}_{pk}(t) = \begin{cases} \sqrt{\frac{c}{2} \frac{\tilde{f}(t)}{\gamma+(I+2)\lambda}}, t \in [t_c - T, t_e - IT - T], \\ \sqrt{\frac{c}{2} \frac{\tilde{f}(t)}{\gamma+(I+1)\lambda}}, t \in (t_e - IT - T, t_c) \end{cases}. \quad (28a)$$

where $\tilde{f}(t)$ denotes *the total boarding demand* at undersaturated stops served by the dispatch at $(t + iT), t \in [t_c - T, t_c), i = 0,1, ...,$ which is given by

$$\tilde{f}(t) = \begin{cases} \sum_{i=0}^{I+1} \sum u_{(t+iT)} f_{u_{(t+iT)}}(t+iT), t \in [t_c - T, t_e - IT - T] \\ \sum_{i=0}^{I} \sum u_{(t+iT)} f_{u_{(t+iT)}}(t+iT), t \in (t_e - IT - T, t_c) \end{cases}. \quad (28b)$$

**Proposition 2**. During the pre-critical period $[t_c - T, t_c)$, there is no queue at any stops along the line.

**Proof**. Proposition 2 will be proved by Lemma 1 in next section. ∎

Per Proposition 2, the optimal peak-period dispatching rate is

$$g_{pk}^*(t) = \begin{cases} \max\left(\sqrt{\frac{c}{2}\frac{\tilde{f}(t)}{\gamma+(I+2)\lambda}}, r_m(t)\right), t \in [t_c - T, t_e - IT - T], \\ \max\left(\sqrt{\frac{c}{2}\frac{\tilde{f}(t)}{\gamma+(I+1)\lambda}}, r_m(t)\right), t \in (t_e - IT - T, t_c) \end{cases}, \quad (29a)$$

where $\tilde{f}(t)$ becomes

$$\tilde{f}(t) = \begin{cases} \sum_{j \in J} f_j(t) + \sum_{i=1}^{I+1} \sum_{u_{(t+iT)}} f_{u_{(t+iT)}}(t+iT), t \in [t_c - T, t_e - IT - T] \\ \sum_{j \in J} f_j(t) + \sum_{i=1}^{I} \sum_{u_{(t+iT)}} f_{u_{(t+iT)}}(t+iT), t \in (t_e - IT - T, t_c) \end{cases}. \quad (29b)$$

Comparing (29a) against (12) finds that the optimal peak-period dispatching rates of general bus lines and shuttle/feeder lines are in a very similar form. The result of the latter is actually a special case of (29) when $r_m(t) = \tilde{f}(t) = f(t)$ for shuttle/feeder lines. Similar to (12), (29) also has a special form of practical importance, i.e.,

$$g_{pk}^*(t) = r_m(t), t \in [t_c - T, t_c), \quad (30)$$

when the peak demand is high and $\gamma, \lambda \gg 1$. However, this result should be used with caution. There is the possibility of $\tilde{f}(t)$ being much larger than $r_m(t)$; and consequently, (30) may not hold anymore.

### 4.1.3. The optimal fleet size

Under the optimal peak-period dispatching profile, $g_{pk}^*(t)$, the solution of $M^*$ can be found in a similar manner to that in Section 3.2. Specifically, we have three unknown variables, i.e., $M, t_c, t_e$. They must satisfy the first-order condition,

$$\frac{dJ_{pk}}{dM} = 0 \Rightarrow$$

$$\gamma - \int_{t_c-T}^{t_c} 0 \left|p_{(t)}\right| dt - \int_{t_c-T}^{t_c} \sum_{i=1}^{I} (i - i_{p,q}) \left|p_{(t+iT)}\right| dt - \int_{t_c-T}^{t_e-(I+1)T} (I + 1 - i_{p,q}) \left|p_{(t+IT+T)}\right| dt = 0 \Rightarrow,$$

$$\gamma = T \sum_{i=1}^{I} (i - i_{p,q}) \cdot \overline{p_{(t+iT)}} + (t_e - t_c - IT)(I + 1 - i_{p,q}) \cdot \overline{p_{(t+IT+T)}}, \quad (31)$$

where $\left|p_{(t+iT)}\right|$ yields the number of the oversaturated stops at time $t + iT, t \in [t_c - T, t_c), 0 \leq i_{p,q} \leq i = 0,1,,\ldots, I$ and at time $t + IT + T, t \in [t_c - T, t_e - (I + 1)T]$; and $\overline{p_{(t+iT)}} \equiv \frac{\int_{t_c-T}^{t_c} |p_{(t+iT)}| dt}{T}$ means the average number of oversaturated stops during $[t_c - T + iT, t_c + iT)$; similarly, $\overline{p_{(t+IT+T)}} \equiv \frac{\int_{t_c-T}^{t_e-IT-T} |p_{(t+IT+T)}| dt}{t_e-t_c-IT}$.

**Lemma 1**. In (31), $|p_{(t)}|$ or $\overline{p_{(t)}}, t \in [t_c - T, t_c)$ (when $i = 0$) does not impact $\frac{dJ_{pk}}{dM}$. It implies that they must satisfy $|p_{(t)}| = \overline{p_{(t)}} = 0$ such that waiting time in queue during $t \in [t_c - T, t_c)$ is minimized with everything else unchanged (no additional cost caused to transit operator and conditions (27a, b) remain being satisfied).

From the definition of $t_c$, we have a boundary condition,

$$G(t_c) - G(t_c - T) = M. \tag{32}$$

Another boundary condition is provided by Proposition 3.

**Proposition 3.** At the end of the critical time period $t_e$, the following boundary condition hold

$$R_M(t_e) \equiv \int_{t_c}^{t_e} \max_{j \in J}(\{r_j(t)\}) \, dt = G(t_e) - G(t_c) \Rightarrow,$$

$$R_M(t_e) = \int_{t_c-T}^{t_e-IT-T} g(t) \, dt + IM = (I+1)M - \int_{t_e-IT-T}^{t_c} g(t) \, dt. \tag{33}$$

which implies that queue(s) at any stops forms at $t_c$ and dissipates at time $t_e$.

**Proof.** Proposition 3 is proved by the following Lemma 2 and 3. ∎

**Lemma 2.** During critical period $[t_c, t_e]$, there exists no time interval with vacant seats departing from the critical stop(s) with peak cross-sectional demand flux along the line.

**Proof.** Without loss of generality, assume there are $\Delta M$ vacant seats observed at the critical stop with peak flux during a certain interval $\Delta t \in [t_c, t_e]$. Then, $\Delta M$ seats can always be dispatched earlier or later to shorten queues. Doing so will only reduce the waiting time in queues with other costs unchanged in $\mathcal{J}_{pk}$. Thus, $\Delta M$ has to be zero. ∎

**Lemma 3.** At time $t_e$, there exists no queue(s) at any stops along the line.

**Proof.** Without loss of generality, assume at time $t_e$ there still exists a queue that lasts for a short interval $\Delta t > 0$ beyond $t_e$. Since beyond $t_e$ is off-peak period when the fleet size is not binding, i.e., $G(t_e + \Delta t) < G(t_e + \Delta t - T) + M$. Per Assumption 1, if any queues remain, a portion of idling fleet ($G(t_e + \Delta t - T) + M - G(t_e + \Delta t) > 0$) can be spared to dispatch to further shorten the queue during $[t_e, t_e + \Delta t]$ until $\Delta t = 0$ and there is no queue at $t_e$. ∎

Substituting $g_{pk}^*(t)$ of (29) to solve three equations (31-33) will produce the optimal fleet size, $M^*$, and the critical time period, $[t_c, t_e]$.

### 4.1.4. Solution method

We note in (29) that $g_{pk}^*(t)$ and $f_{u_{(t+iT)}}(t + iT)$ are interdependent, and thus propose the following iteration method to iteratively find both $g_{pk}^*(t), t \in [t_c - T, t_c)$ and $M^*$.

*Initialization*: Initialize $g^{(0)}(t)$, e.g., $g^{(0)}(t) = \sum_{j=1}^{J} f_j(t)$, and $M^{(0)}$. Set iteration variable $n = 1$.

*Step 1*: Determine $t_c^{(n)}$ by (33) with $g^{(n-1)}(t)$ and $M^{(n-1)}$.

$$t_c^{(n)} | \int_{t_c-T}^{t_c} g^{(n-1)}(t) dt = M^{(n-1)}. \tag{34}$$

*Step 2*: Determine $t_e^{(n)}, I^{(n)}$ by (35) with $g^{(n-1)}(t)$ and the updated $t_c^{(n)}$.

$$t_e^{(n)}, I^{(n)}|\gamma = T \sum_{i=1}^{I}(i - i_{p,q}) \cdot \overline{p_{(t+iT)}} + \left(t_e - t_c^{(n)} - IT\right)(I + 1 - i_{p,q}) \cdot \overline{p_{(t+IT+T)}}, \quad (35)$$

where $\overline{p_{(t+iT)}}, t \in \left[t_c^{(n)} - T, t_c^{(n)}\right), 0 \leq i_{p,q} \leq i = 1,2,...,$ are determined by $g^{(n-1)}(t)$.

*Step 3*: Determine $\tilde{f}(t), t \in [t_c - T, t_c)$ by (36) with $g^{(n-1)}(t)$ and the updated $t_c^{(n)}, t_e^{(n)}, I^{(n)}$.

$$\tilde{f}(t) = \tilde{f}(t) + \delta_{j,u} f_{j_{(t+iT)}}(t + iT), j = 1,...,J, i = 1,2,...I, \quad (36)$$

where $\delta_{j,u} = 1$ if $t_{j,q}$ does not exist yet, or $\int_{t_{j,q}}^{t}\left(r_j(s) - g(s)\right)ds \leq 0$ meaning the queue had dissipated; otherwise, $\delta_{j,u} = 0$.

*Step 4*: Update $g^{(n)}(t)$ according to (29).

*Step 5*: Determine an auxiliary fleet size $\tilde{M}$ by (37) with the updated $g^{(n)}(t), t_c^{(n)}$, and $t_e^{(n)}, I^{(n)}$.

$$\tilde{M} = \frac{R_M\left(t_e^{(n)}\right) + \int_{t_e^{(n)} - IT - T}^{t_c^{(n)}} g^{(n)}(t)dt}{(I+1)}. \quad (37)$$

*Step 6*: Update $M^{(n)}$ by the method of successive averages (MSA)*,

$$M^{(n)} = M^{(n-1)} + \frac{\tilde{M} - M^{(n-1)}}{n}. \quad (38)$$

*Step 6*: Stop if $\frac{\int_{t_c-T}^{t_c}|g^{(n)}(t) - g^{(n-1)}(t)|dt}{M^{(n)}} \leq \varepsilon$, where $\varepsilon$ is a predefined small value, e.g., $\varepsilon = 10^{-3}$; otherwise, let $n = n + 1$ and go to Step 1.

*Note: MSA is employed here to stabilize the convergence process; otherwise, oscillation may occur since $g^{(n)}(t)$ and $M^{(n)}$ are updated in a sequential manner.

### 4.1.5. Numerical examples

A simple transit line with three stops, as shown in Figure 4, is used to demonstrate the application of the above model and solution algorithm. The total OD demand between stops, as indicated by arrow lines in Figure 4, is set as $\{D_{12}, D_{13}, D_{23}\} = \{600, 200, 200\}$ (patrons), where subscripts indicate the origin and destination stop pairs. As the baseline, they all follow the temporal distribution pattern defined by (19) with the same parameter values as given in Table 2. These baseline values will be varied in the following tests.

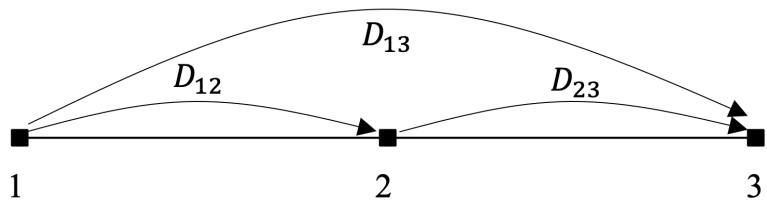

**Figure 4. A transit line with multiple OD.**

First, we show the convergence treads of the solution, $M^*$, in Figure 5. As seen, the convergence is soon reached in 10 iterations with different initial points.

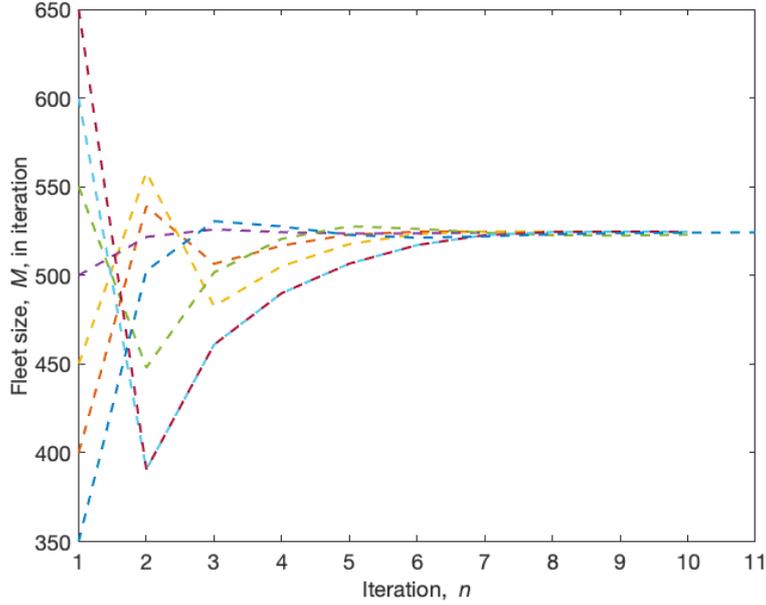

**Figure 5. Convergence of the optimized fleet size with different initial points.**

Figures 6a and 6b visualize the optimized dispatching profiles against the demand profiles at stop 1 and 2, respectively. It is found that queues form at both stops at time $t_c = 77.3$ mins. The queue at stop 1 lasts 25.4 mins and dissipates at time $t_e = 102.7$ mins. The queue at stop 2 lasts 4.6 mins and dissipates at time $t_{e2} = 81.9$ mins. Adding up the queueing durations at two stops yields $25.4 + 4.6 = 30 = \gamma$ as dictated by (31). For the system, however, the queueing duration is $\mathcal{T}_q = 25.4$ mins. Observation verifies that $g_{pk}^*(t) = r_m(t) = r_1(t)$ for $t \in [t_c - T, t_c)$. In this case, the optimized fleet size is $M^* = 525$ seats.

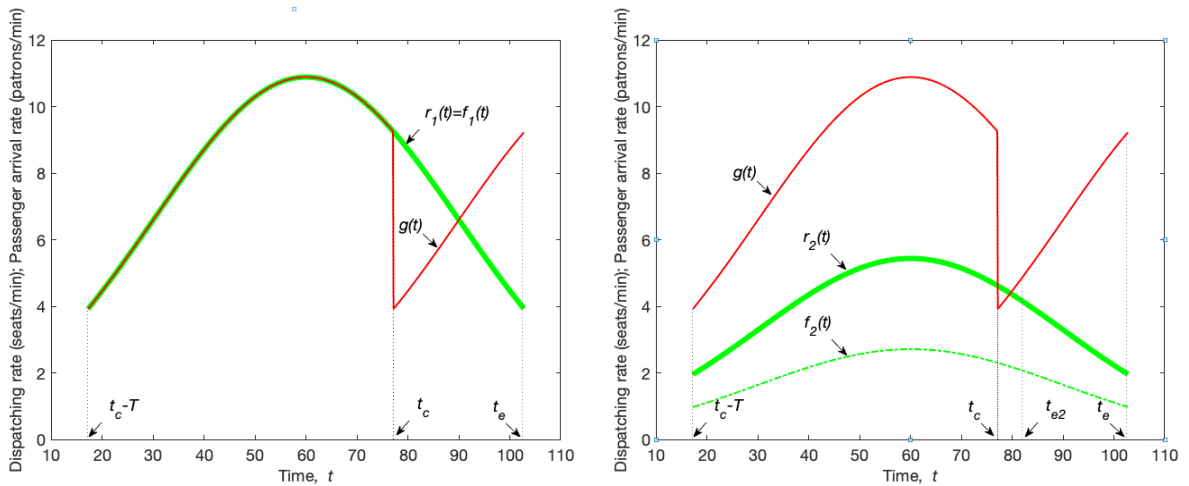

a. Dispatching and demand profiles at stop 1   b. Dispatching and demand profiles at stop 2

**Figure 6. Results with the baseline setting.**

We also test scenarios with two stops' demand peaking at different times. Specifically, we set $\mu_{23} \in \{20,100\}$ (mins) to represent that the demand at stop 2 peaks at 20 and 100 mins, respectively. With other parameters remaining the same with the baseline settings, the optimized results are shown in Figures 7a-d. As seen in Figures 7a and 7b for $\mu_{23} = 20$ mins, no queue forms at stop 2 due to the shift of the demand peak time. As a result, the optimized dispatching profile is determined by $g_{pk}^*(t) = r_m(t) = \max(r_1(t), r_2(t)), t \in [t_c - T, t_c)$, which explains the beginning branch of $g_{pk}^*(t)$ profile for a short interval after $t_c - T$. In this case, the queueing period begins earlier at 74.7 mins and ends at 104.7 mins with a duration of 30 mins. The optimized fleet size is reduced to $M^* = 513$ seats.

Figures 7c and 7d are for the case of $\mu_{23} = 100$ mins. Comparatively, the results differ little from that in Figure 6. Both stop 1 and 2 witness queues that begin at $t_c = 77.0$ mins but last to $t_e = 103.0$ and $t_{e2} = 81.0$ mins, respectively. The optimized fleet size becomes $M^* = 523$ seats. In addition, it is found that both cases of $\mu_{23} \in \{20,100\}$ lengthen the queueing duration of the system, with a different extent though. This is not surprising since the dispersing of demand peaks at stops tends to stretch the saturation state of the system.

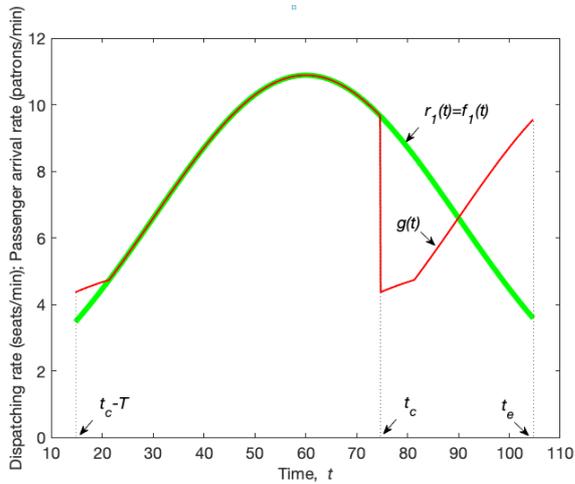

a. Dispatching and demand profiles at stop 1, $\mu_{23} = 20$ (mins)

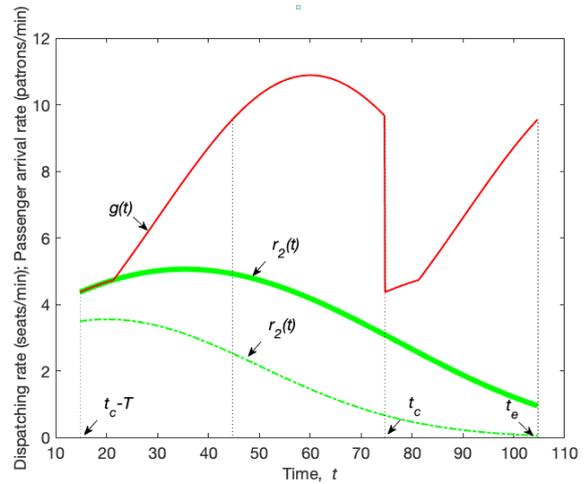

b. Dispatching and demand profiles at stop 2, $\mu_{23} = 20$ (mins)

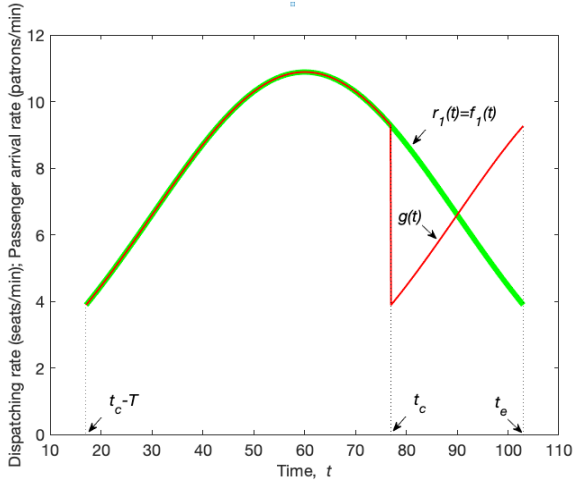
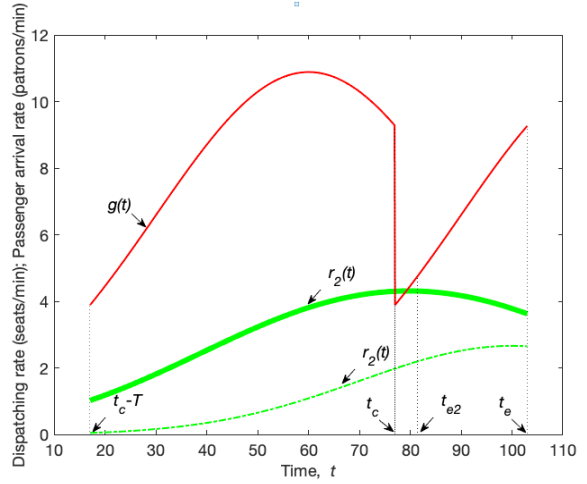

c. Dispatching and demand profiles at stop 1, $\mu_{23} = 100$ (mins)

d. Dispatching and demand profiles at stop 2, $\mu_{23} = 100$ (mins)

**Figure 7. Results with different demand peak times at stop 2.**

Lastly, we examine a costly-fleet scenario with $\gamma = 120$ (mins/seat) and $\mu_{12}, \mu_{13}, \mu_{23} = 70$ (mins) (for keeping $t_c - T \geq 0$). Figures 8a and 8b present the dispatching profiles at two stops, respectively. In contrast to the results reported in Figure 6, the critical-period duration in this scenario lasts longer than one $T$ and the beginning time ($t_c = 68.0$ mins) starts earlier (even with the demand peak time postponed from 60 to 70 mins). The queueing duration at stop 1 and 2 are 81 and 39 mins, respectively, which again add up to $120 = \gamma$. In this case, the optimized fleet size is $M^* = 368$ seats.

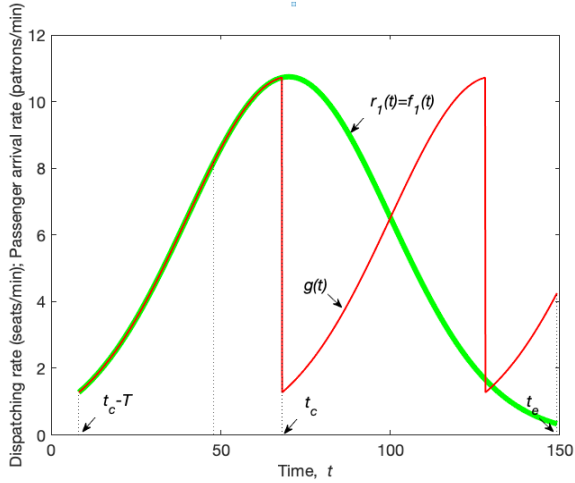
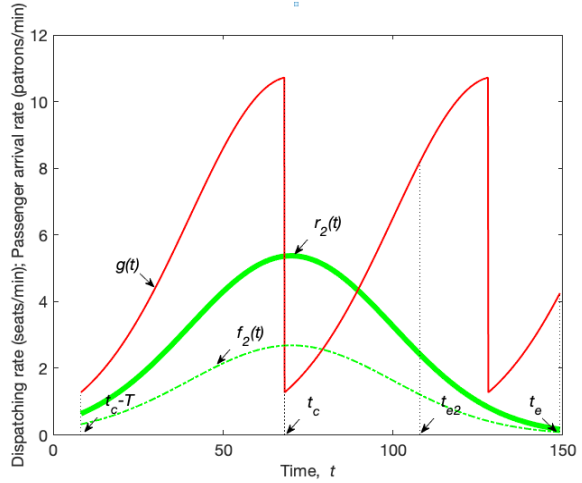

a. Dispatching and demand profiles at stop 1   b. Dispatching and demand profiles at stop 2

**Figure 8. Results with costly fleet, $\gamma = 120$ (mins/seat).**

Again, we find that the divergence between the special result (30) and the general form (29) of the optimal peak-period dispatching rate only occurs at very low demand levels. The detailed results are omitted.

## 4.2. Scheduling with mixed-size or modular buses

This section demonstrates the scheduling for mixed-size or modular buses using Hurdle's original model (1). Similar results can also be obtained for the case of general bus lines, but omitted here for the sake of brevity.

We note that the operator's cost parameters, i.e., $\gamma$ and $\lambda$, may not be constant but rather depends on the size of vehicles, i.e., $c$ (seats/vehicle). On realizing the economies of scale in transit industry, researchers had introduced functions describing the varying costs related to vehicle sizes. For instance, a commonly used one is in the following form that is adapted to the unit of cost per seat (Chang and Schonfeld, 1991; Kim and Schonfeld, 2019):

$$\gamma = \frac{\gamma_0(c_0+\alpha c)}{c} \text{ and } \lambda = \frac{\lambda_0(c_0+\alpha c)}{c}, \tag{39}$$

where $c_0, \gamma_0, \lambda_0$, and $\alpha$ are known parameters.

Substituting (39) in (18) for the off-peak period yields the optimization problem with respect to $c$,

$$\begin{aligned}\text{minimize}_c \; \mathcal{J}_{op} &= \int_{t\in \mathbf{E}_{op}} \left( \frac{\lambda_0(c_0+\alpha c)}{c} g_{op}^*(t) + \frac{cf(t)}{2g_{op}^*(t)} \right) dt \\ &= \int_{t\in \mathbf{E}_{op}} \sqrt{2\lambda_0(c_0+\alpha c)f(t)}\, dt,\end{aligned} \tag{40a}$$

subject to:

$$c \geq 2\frac{\lambda_0(c_0+\alpha c)}{c} f(t), t \in \mathbf{E}_{op}, \tag{40b}$$

$$\int_{t-T}^{t} \sqrt{\frac{cf(t)}{2\frac{\lambda_0(c_0+\alpha c)}{c}}} \leq M, t \in \mathbf{E}_{op}, \tag{40c}$$

where constraint (40b) guarantees no queue forms during the off-peak period, to which verification can be done by observing $g_{op}^*(t) = \sqrt{\frac{cf(t)}{2\frac{\lambda_0(c_0+\alpha c)}{c}}} \geq f(t)$; and (40c) is the fleet size constraint.

Combining the first-order condition of (40a) and the boundary condition (40b) yields the solution,

$$\tilde{c}_{op}(t) = \lambda_0 \alpha f(t) + \sqrt{2\lambda_0 c_0 f(t) + \big(\lambda_0 \alpha f(t)\big)^2}, t \in \mathbf{E}_{op}. \tag{41}$$

With $\tilde{c}_{op}$ in (41), the off-peak dispatch rate $g_{op}^*(t)$ will always be $f(t)$ and (40c) is always satisfied (since $f(t)$ in off-peak period is no larger than that in peak period).

Further considering practical reasons, the vehicle size cannot be continuously varying with time but probably would be limited in a finite set, e.g., $\mathbf{c} \equiv \{1\Delta c, 2\Delta c, \ldots, \mathcal{M}\Delta c\}$, where $\Delta c$ can be seen as the size of a modular vehicle and $\mathcal{M}\Delta c$ is the vehicle-size upper bound. Ultimately, we have the optimal off-peak vehicle size as

$$c_{op}^*(t) = \min\left(\max\left(\left\lceil\frac{\tilde{c}_{op}(t)}{\Delta c}\right\rceil, 1\right), \mathcal{M}\right)\Delta c, t \in \mathbf{E}_{op}, \qquad (42)$$

where $\lceil \cdot \rceil$ returns the closest integer no less than the argument, and the meaning is the group size of modular vehicles.

For the peak period $[t_c - T, t_e]$, assuming that $g_{pk}^*(t) = f(t), t \in [t_c - T, t_c)$ holds, we have the following optimization problem with respect to $c$,

$$\begin{aligned}\text{minimize}_c \; \mathcal{J}_{pk} &= \frac{\gamma_0(c_0+\alpha c)}{c}M + \frac{\lambda_0(c_0+\alpha c)}{c}\int_{t_q-T}^{t_q+\mathcal{T}_q} g_{pk}^*(t)dt + \int_{t_q-T}^{t_q} \frac{cf(t)}{2g_{pk}^*(t)}dt \\ &\quad + \int_{t_q}^{t_q+\mathcal{T}_q}\left(\frac{c}{2} + F(t) - G_{pk}^*(t)\right)dt \\ &= \frac{\gamma_0(c_0+\alpha c)}{c}\int_{t_q-T}^{t_q} f(t)dt + \frac{\lambda_0(c_0+\alpha c)}{c}\Big[(I+1)\int_{t_q-T}^{t_q} f(t)dt + \\ &\quad \int_{t_q-T}^{t_q+\mathcal{T}_q-IT-T} f(t)dt\Big] + \Big[\frac{(T+\mathcal{T}_q)c}{2} + \int_{t_q}^{t_q+\mathcal{T}_q}\left(F(t)-G_{pk}^*(t)\right)dt\Big].\end{aligned} \qquad (43)$$

The first-order condition to (43), $\frac{d\mathcal{J}_{pk}}{dc} = 0$, yields,

$$\tilde{c}_{pk}(t) = \begin{cases} \sqrt{\frac{2(\gamma_0 c_0 + \lambda_0 c_0(I+2))f(t)}{I+2}}, & t \in [t_q - T, t_q + \mathcal{T}_q - IT - T] \\ \sqrt{\frac{2(\gamma_0 c_0 + \lambda_0 c_0(I+1))f(t)}{I+1}}, & t \in (t_q + \mathcal{T}_q - IT - T, t_q) \end{cases}. \qquad (44)$$

Again, combining practical constraints regarding modular vehicles, we have the optimal peak-period vehicle size as,

$$c_{pk}^*(t) = \min\left(\max\left(\left\lceil\frac{\tilde{c}_{pk}(t)}{\Delta c}\right\rceil, 1\right), \mathcal{M}\right)\Delta c, t \in [t_c - T, t_c), \qquad (45a)$$

$$c_{pk}^*(t) = c_{pk}^*(t - T), t \in [t_c, t_e]. \qquad (45b)$$

Observing (41, 44) reveals some insights: (i) The variant-cost parameter, $\alpha$, affects only the optimal off-peak vehicle size, but has nothing to do with the peak-period result; and (ii) larger-sized buses will be preferred if purchasing and operating buses are relatively more expensive, i.e., $\gamma_0, \lambda_0, c_0$ have large values.

We note the above results are dictated by the operator's cost functions given in (39); and alternative cost functions exist in the literature, e.g., $\gamma = \gamma_0(c_0 + c^\alpha)/c$ and $\lambda = \lambda_0(c_0 + c^\alpha)/c$ where $\alpha < 1$ (Chen et al., 2020). To different applications, the above solution procedure still applies. Although closed-form solutions will not be guaranteed, numerical solutions could be readily found by solving the first-order conditions (as higher-degree polynomial equations) of the scheduling problem, thanks to the proposed mathematical formulation.

## 5. Conclusions

This paper completes the scheduling problem of Hurdle (1973) by proposing a novel mathematical solution approach. The original constrained variational problem is relaxed to an unconstrained one and solved using calculus of variations. Hurdle's finding of the optimal peak-period dispatching rate is found to be a special result of the general form obtained by the Euler-Lagrange equation of the relaxed problem. The unsolved optimal fleet size is now solved.

Based on the proposed approach, we make two extensions to the scheduling problem of general bus lines with multiple origins and destinations and that of mixed-size or modular buses. Closed-form solutions are derived with new insights uncovered. Among others, we show that the optimal pre-queueing/critical dispatching rate for a general line is oftentimes determined by the maximum cross-sectional demand flux along the line.

Numerical examples demonstrate the effectiveness of the proposed approach. Parameter analyses are also conducted thanks to the high computation efficiency of our solution method. A useful by-product is the queueing profile under the optimal schedule design. Accordingly, proper measures can be designed to manage/guide the queue.

Of note, the above modeling of the schedule problems still has several limitations. For instance, the current service of interest is still limited to a single line. The pre-given and fixed cycle time may vary in real world for different times of the day, depend on the number of boarding and alighting patrons, and even be a stochastic variable due to random influential factors (e.g., signals and drivers' behaviors). Last but not least, patrons' behaviors are much simplified in the current study, which would not be well justified in the cases of common lines, alternative travel modes, heterogeneous patrons, etc. Some of these limitations promise extensions. Select topics are under exploration.

**Acknowledgments**

The research was supported by a fund provided by the National Natural Science Foundation of China (No. 51608455) and Sichuan Science & Technology Program (No. 2020YFH0038).